\documentclass[reqno]{amsart}

\theoremstyle{theorem}
\newtheorem{theorem}{Theorem}[section]
\newtheorem{corollary}[theorem]{Corollary}
\newtheorem{lemma}[theorem]{Lemma}
\newtheorem{prop}[theorem]{Proposition}

\theoremstyle{definition}
\newtheorem{definition}{Definition}[section]
\newtheorem{remark}[definition]{Remark}
\newtheorem{example}[definition]{Example}

\numberwithin{equation}{section}

\newcommand{\PP}{\mathbb{P}}
\newcommand{\CC}{\mathbb{C}}
\newcommand{\RR}{\mathbb{R}}
\newcommand{\QQ}{\mathbb{Q}}
\newcommand{\ZZ}{\mathbb{Z}}
\newcommand{\BA}{\mathbb{A}}
\newcommand{\KK}{\mathbb{K}}

\newcommand{\CP}{\mathcal{P}}
\newcommand{\CE}{\mathcal{E}}
\newcommand{\CF}{\mathcal{F}}
\newcommand{\CO}{\mathcal{O}}
\newcommand{\CM}{\mathcal{M}}
\newcommand{\CA}{\mathcal{A}}
\newcommand{\CL}{\mathcal{L}}

\newcommand{\CJ}{\mathcal{J}}

\newcommand{\PD}[1]{\frac{\partial}{\partial {#1}}}
\newcommand{\D}{\delta}

\begin{document}
\title[Floer cohomology of Lagrangian torus fibers]
{Floer cohomology and disc instantons\\
of Lagrangian torus fibers in\\
Fano toric manifolds}

\author{Cheol-Hyun Cho
\and Yong-Geun Oh}
\thanks{The second named author is partially supported
by the NSF grant \#DMS 0203593, a grant of the 2000 Korean Young
Scientist Prize, and the Vilas Research Award of University of
Wisconsin}

\address{Department of Mathematics, Northwestern University,
Evanston, IL 60208, cho@math.wisc.edu\\
Department of Mathematics, University of Wisconsin, Madison, WI
53706 \& Korea Institute for Advanced Study, 207-43
Cheongryangri-dong Dongdaemun-gu, Seoul 130-012, KOREA,
oh@math.wisc.edu}

\begin{abstract}
In this paper, we first provide an explicit description of {\it
all} holomorphic discs (``disc instantons'') attached to
Lagrangian torus fibers of arbitrary compact toric manifolds, and
prove their Fredholm regularity. Using this, we compute
Fukaya-Oh-Ohta-Ono's (FOOO's) obstruction (co)chains and the Floer
cohomology of Lagrangian torus fibers of Fano toric manifolds. In
particular specializing to the formal parameter $T^{2\pi} =
e^{-1}$, our computation verifies the folklore that FOOO's
obstruction (co)chains correspond to the Landau-Ginzburg
superpotentials under the mirror symmetry correspondence, and also
proves the prediction made by K. Hori about the Floer cohomology
of Lagrangian torus fibers of Fano toric manifolds. The latter
states that the Floer cohomology (for the parameter value
$T^{2\pi} = e^{-1}$) of all the fibers vanish except at a finite
number, the Euler characteristic of the toric manifold, of base
points in the momentum polytope that are critical points of the
superpotential of the Landau-Ginzburg mirror to the toric
manifold. In the latter cases, we also prove that the Floer
cohomology of the corresponding fiber is isomorphic to its
singular cohomology.

We also introduce a restricted version of the Floer cohomology of
Lagrangian submanifolds, which is a priori more flexible to define
in general, and which we call the {\it adapted Floer cohomology}.
We then prove that the adapted Floer cohomology of any
non-singular torus fiber of Fano toric manifolds is well-defined,
invariant under the Hamiltonian isotopy and isomorphic to the
Bott-Morse Floer cohomology of the fiber.
\end{abstract}

\maketitle

\bigskip

\section{Introduction}

Floer cohomology of Lagrangian intersections was introduced by
Floer \cite{Fl} in symplectic geometry. Since then, its
construction has been further generalized \cite{Oh1} and an
obstruction theory to its definition has been developed by
Fukaya-Oh-Ohta-Ono \cite{FOOO}. It has been proven to be a
powerful tool in studying various problems in symplectic geometry
(see \cite{Fl}, \cite{Oh4}, \cite{Che}, \cite{P}, \cite{Se},
\cite{FOOO}, \cite{BC}, and \cite{TY}, for example). The theory
itself was greatly enhanced by the advent of the Fukaya category
\cite{Fuk1} and the homological mirror symmetry proposal by
Kontsevich \cite{Ko}, and also by the open string theory of
$D$-branes in many physics literature, among which
\cite{hori:ms02}, \cite{hori:lmsd00} will be the most relevant to
the content of the present paper.

Even in the midst of these theoretical enhancement and successful
applications of the Floer theory, actual computation of Floer
cohomology itself for specific examples remains to be a
non-trivial task, especially with $\ZZ$-coefficients (not just
with $\ZZ_2$-coefficients), except for the cases where there is no
quantum contribution \cite{Fl} or for the case of {\it real
manifolds} i.e., the fixed point sets of anti-holomorphic
involutions \cite{Oh2}, \cite{FOOO}. Indeed, computation of the
Floer cohomology in the presence of nontrivial holomorphic discs
requires detailed understanding of the quantum contribution of the
holomorphic discs (or the effect of ``open string instantons'' in
the physics terminology) to the cohomology of the Lagrangian
submanifolds. In this respect, the recent computation
\cite{cho:hsf03} by the first named author of the Floer cohomology
of the Clifford torus in $\PP^n$ sheds some light on a general
procedure of computing the Floer cohomology ``by direct
calculation of disc instanton effects'' in the context of
$A$-model without relying on the $B$-model calculations and the
mirror symmetry correspondence, which is still conjectural.

In this paper, we extend this computation and compute the ({\it
adapted}) Floer cohomology of {\it all} the non-singular torus
fibers of smooth Fano toric varieties equipped with symplectic
(K\"ahler) form. As in \cite{cho:hsf03}, we will carry out this by
computing the Bott-Morse version $HF^{BM}(L;J_0)$ of the Floer
cohomology of Lagrangian submanifold $L$ that was introduced in
\cite{FOOO}. Our computation, when the Floer cohomology is twisted
with the flat line bundles and the formal parameter $T$ is set
$T^{2\pi} = e^{-1}$, verifies the prediction made by Hori-Vafa
\cite{hori:ms02} for the Lagrangian torus fibers of Fano toric
manifolds based on the mirror symmetry correspondence via the
linear sigma models \cite{Wi}, \cite{hori:ms02}.

In the point of view of the obstruction theory developed in
\cite{Oh1}, \cite{FOOO}, a priori, the torus fibers of general
toric manifolds are neither {\it monotone} nor {\it unobstructed},
and may carry holomorphic discs of non-positive Maslov indices.
Recall that the Clifford torus is also obstructed {\it as an
object} in the $A_\infty$-category \cite{Oh1}, \cite{FOOO}, but
the fact that it is monotone enables one to define the Floer
cohomology \cite{Oh1},\cite{Oh4} which the first named author
computed in \cite{cho:hsf03}. Combination of these facts prevent
us from directly applying the general construction of the Floer
cohomology from \cite{FOOO} and forces us to manually construct a
restricted version of the Floer cohomology and to prove the
invariance property. For this purpose, some specific geometry of
the moduli of holomorphic discs associated to the pair $(L,J_0)$
of the torus fiber $L$ and the canonical complex structure $J_0$
on the toric variety will play an essential role both for the
definition and computation of the Floer cohomology. We will prove
that there exists no non-constant holomorphic discs of
non-positive Maslov indices for the torus fibers, although its
Hamiltonian deformations of them may allow such
(pseudo-)holomorphic discs. Our definition of the adapted Floer
cohomology exploits this specific feature of the pair $(L,J_0)$.
We call this version of the Floer cohomology the {\it adapted
Floer cohomology}. It appears that in general this adapted Floer
cohomology is more flexible to define and exploits best specific
features of the moduli of holomorphic discs of the given pair
$(L,J_0)$. (However the arguments from \cite{FOOO} involving the
homotopy inverse of the $A_\infty$-algebra strongly suggests that
whenever the adapted Floer cohomology is defined, the deformed
Floer cohomology in the sense of \cite{FOOO} will be also
well-defined and isomorphic to the adapted Floer cohomology. This
question will be studied in the final version of \cite{FOOO}.)

Once the well-definedness of the adapted Floer cohomology is
established, its computation largely follows the scheme used by
the first named author \cite{cho:hsf03}: Firstly, we derive
general Maslov index formula of holomorphic discs in terms of the
intersection number of natural divisors associated to the toric
manifolds. Secondly we explicitly classify all the holomorphic
discs and prove the Fredholm regularity of the discs. Then using
this information, we compute the Bott-Morse version
$HF^{BM}(L;J_0)$ of the Floer cohomology of $L$ with respect to
the complex structure $J_0$. Because the torus fibers do not have
non-constant holomorphic discs of non-positive Maslov indices (see
Theorem \ref{classify}), the argument from \cite{Oh4} proves
$HF^{BM}(L;J_0)$ isomorphic to the adapted Floer cohomology
$HF^{ad}(L;J_0)$ in the Fano case, and expected to be the same in
general.

In the course of our computation, we also derive an area formula
Theorem \ref{area} for the holomorphic discs of the Maslov index 2
(and so of {\it all} holomorphic discs) in terms of the location
of the base point of the Lagrangian fiber and the relative
homology class of the disc (or the divisor of the toric manifold
that the disc intersect). This formula is crucial for our proof of
the prediction that the base points in the momentum polytope at
which the corresponding fiber has non-trivial Floer cohomology are
indeed those corresponding to the critical points of the
superpotential of the Landau-Ginzburg mirror.

We would like to emphasize that the mirror symmetry prediction
made both in the Kontsevich proposal or by physicists does not
really concern the standard Floer cohomology in symplectic
geometry which uses the Novikov ring as its coefficients, but its
convergent power series version. One byproduct of our
classification of disc instantons is that this latter version of
the Floer cohomology is defined and so substitution of the formal
parameter $T^{2\pi}$ by the number $e^{-1}$ is allowed in the Fano
toric  case. However the latter version of the Floer cohomology is
{\it not} known to be invariant in general under the Hamiltonian
isotopy of the Lagrangian torus fiber and so the mirror symmetry
prediction concerns the {\it K\"ahler geometry} of the Lagrangian
torus fibers (with respect to the natural complex structure $J_0$
and the K\"ahler form $\omega$), rather than the {\it symplectic
geometry} of its Hamiltonian isotopy class. For example, it is
possible that a fiber has trivial Floer cohomology with Novikov
ring as its coefficients but non-trivial one with the parameter
value $T^{2\pi} = e^{-1}$ (see section 13 for an explicit example
of Hirzebruch surfaces).

Our work also provides some concrete mathematical evidence in the
toric case for the conjectural relation between the superpotential
and the ``open Gromov-Witten invariants'' which has been advocated
by physicists (see \cite{KKLM} for example). More precisely, we
verifies that under the mirror symmetry correspondence of a torus
fiber, the {\it one-point} open Gromov-Witten invariant, which is
essentially FOOO's obstruction chain [FOOO], maps to the
superpotential $W$ of the Landau-Ginzburg mirror, and {\it
two-point} invariants, which is essentially the Floer differential
$\delta_2\langle pt \rangle$ in the Bott-Morse setting, maps to
the derivative $\frac{\partial W}{\partial \Theta}$. We refer to
section \ref{sec:obstruct} for more discussion on this point.

One general distinction between the Fano and the non-Fano cases
lies in the transversality property of the singular strata of
various compactified moduli spaces. More precisely, non-Fano
manifolds carry spheres of {\it negative} Chern numbers and so the
compactified moduli space may contain singular strata that contain
sphere bubbles (especially their multiple covers) of negative
Chern numbers. As the study in \cite{FOOO} demonstrated, such
problems in the moduli space of holomorphic discs in relation to
the Floer theory (or to open Gromov-Witten invariants) are much
more troublesome than the case of spheres. We refer to section
\ref{sec:discussion} for more detailed discussion on this.

We like to thank K. Hori for explaining us the mirror symmetry
correspondence via the Landau-Ginzburg models and his $B$-model
calculation that leads to his conjectural description of the Floer
cohomology of the fibers of Fano toric manifolds. The first named
author would like to thank L. Borisov and S. Hu for helpful
discussions. The second named author thanks K. Fukaya, H. Ohta and
K. Ono for their interest in the results of this paper and some
interesting discussion on the homotopy inverse defined in
\cite{FOOO} during his visit of Kyoto University in August, 2003.

\section{Compact toric manifolds}
\label{sec:toric}

We consider smooth and compact toric varieties. Here we closely
follow the Batyrev \cite{batyrev:qcrtm92} with minor notational
changes (See M. Audin \cite{audin:ttasm91} for more details)

In order to obtain an n-dimensional compact toric manifold $V$, we
need a combinatorial object $\Sigma$, a {\em complete fan of
regular cones}, in a $n$-dimensional vector space over $\RR$.

Let $N$ be the lattice $\ZZ^n$, and let $M=Hom_\ZZ(N,\ZZ)$ be the
dual lattices of rank $N$. Let $N_{\RR} = N \otimes \RR$ and
$M_{\RR} = M \otimes \RR$.
\begin{definition}
A convex subset $\sigma \subset N_{\RR}$ is called a regular
$k$-dimensional cone $(k\geq 1)$ if there exists $k$ linearly
independent elements $v_1,\cdots,v_k \in N$ such that
$$ \sigma = \{a_1v_1 + \cdots + a_k v_k \mid a_i \in \RR, a_i \geq
0\},$$ and the set $\{v_1,\cdots,v_k\}$ is a subset of some
$\ZZ$-basis of $N$. In this case, we call $v_1,\cdots,v_k \in N$
the integral generators of $\sigma$.
\end{definition}
\begin{definition}
A regular cone $\sigma'$ is called a {\em face} of  a regular cone
$\sigma$ (we write $\sigma' \prec \sigma$) if the set of integral
generators of $\sigma'$ is a subset of the set of integral
generators of $\sigma$.
\end{definition}
\begin{definition}
A finite system $\Sigma = {\sigma_1,\cdots,\sigma_s}$ of regular
cones in $N_{\RR}$ is called a {\em complete $n$-dimensional fan}
 of regular cones, if the following conditions are satisfied.
\begin{enumerate}
\item if $\sigma \in \Sigma$ and $\sigma' \prec \sigma$, then $\sigma'
\in \Sigma$;

\item if $\sigma, \sigma'$ are in $\Sigma$, then $\sigma' \cap
\sigma \prec \sigma$ and $\sigma' \cap \sigma \prec \sigma'$;

\item $N_{\RR} = \sigma_1 \cup \cdots \cup \sigma_s$.
\end{enumerate}
\end{definition}
The set of all $k$-dimensional cones in $\Sigma$ will be denoted
by $\Sigma^{(k)}.$

\begin{example}\label{ex:cpn}
Consider basis vectors $e_1,\cdots,e_n$ in a $n$-dimensional real
vector space. Let $v_i =e_i$ for $i=1,\cdots,n$ and let $v_{n+1}
=-e_1-e_2-\cdots-e_n$. Any $k$-element subset $I \subset
\{v_1,\cdots,v_{n+1}\}$ for $(k \leq n)$ generates a
$k$-dimensional regular cone $\sigma(I)$. The set $\Sigma(n)$
consisting of $2^{n+1}-1$ cones $\sigma(I)$ generated by $I$ is a
complete $n$-dimensional fan of regular cones, with which later we
will associate a projective space $\PP^n$.
\end{example}

\begin{definition}\label{prim}
Let $\Sigma$ be a complete $n$-dimensional fan of regular cones.
Denote by $G(\Sigma) = \{v_1,\cdots,v_N\}$ the set of all
generators of 1-dimensional cones in $\Sigma$ ( $N=$ Card
$\Sigma^{(1)}$). We call a subset $\CP =
\{v_{i_1},\cdots,v_{i_p}\} \subset G(\Sigma)$ a {\bf primitive
collection} if $\{v_{i_1},\cdots,v_{i_p}\}$ does not generate
$p$-dimensional cone in $\Sigma$, while for all $k \, (0 \leq k <
p)$ each $k$-element subset of $\CP$ generates a $k$-dimensional
cone in $\Sigma$.
\end{definition}
\begin{example}\label{cpnprim}
Let $\Sigma$ be a fan from Example \ref{ex:cpn}. Then there exists
the unique primitive collection $\CP$ which is the set of all
generators $\{v_1,\cdots,v_{n+1}\}$.
\end{example}
\begin{definition}
Let $\CC^N$ be $N$-dimensional affine space over $\CC$ with the
set of coordinates $z_1,\cdots,z_N$ which are in the one-to-one
correspondence $z_i \leftrightarrow v_i$ with elements of
$G(\Sigma)$. Let $\CP = \{v_{i_1},\cdots,v_{i_p}\}$ be a primitive
collection in $G(\Sigma)$. Denote by $\BA(\CP)$ the
$(N-p)$-dimensional affine subspace in $\CC^n$ defined by the
equations
$$z_{i_1}= \cdots=z_{i_p}=0.$$
\end{definition}

\begin{remark}
Since every primitive collection $\CP$ has at least two elements,
the codimension of $\BA(\CP)$ is at least 2.
\end{remark}

\begin{definition}\label{homo}
Define the closed algebraic subset $Z(\Sigma)$ in $\CC^N$ as follows
$$Z(\Sigma) = \cup_{\CP} \BA(\CP),$$
where $\CP$ runs over all primitive collections in $G(\Sigma)$.
Put
$$U(\Sigma) = \CC^N \setminus Z(\Sigma).$$
\end{definition}

\begin{definition}
Let $\KK$ be the subgroup in $\ZZ^N$ consisting of all lattice vectors
$\lambda = (\lambda_1,\cdots,\lambda_N)$ such that
$$\lambda_1 v_1 + \cdots + \lambda_N v_N =0.$$
\end{definition}
Obviously $\KK$ is isomorphic to $\ZZ^{N-n}$ and we have the exact
sequence:
\begin{equation}\label{kexact}
 0 \to \KK \to \ZZ^N \stackrel{\pi}{\to} \ZZ^n \to 0,
\end{equation}
where
the map $\pi$ sends the basis vectors $e_i$ to $v_i$ for
$i=1,\cdots,N$.

\begin{definition}
Let $\Sigma$ be a complete $n$-dimensional fan of regular cones.
Define $D(\Sigma)$ to be the connected commutative subgroup in
$(\CC^*)^N$ generated by all one-parameter subgroups
$$a_{\lambda} : \CC^* \to (\CC^*)^N,$$
$$ t \mapsto (t^{\lambda_1},\cdots,t^{\lambda_N})$$
where $\lambda = (\lambda_1, \cdots, \lambda_N) \in \KK$.
\end{definition}

It is easy to see from the definition that $D(\Sigma)$ acts freely
on $U(\Sigma)$. Now we are ready to give a definition of the
compact toric manifold $X_{\sigma}$ associated with a complete
n-dimensional fan of regular cones $\Sigma$.

\begin{definition}
Let $\Sigma$ be a complete $n$-dimensional fan of regular cones.
Then the quotient
$$X_{\Sigma} = U(\Sigma)/D(\Sigma)$$
is called the {\em compact toric manifold associated with $\Sigma$}.
\end{definition}

\begin{example}
Let $\Sigma$ be a fan $\Sigma(n)$ from Example \ref{ex:cpn}.
By \ref{cpnprim}, $U(\Sigma(n)) = \CC^{n+1} \setminus \{0\}$.
By the definition of $\Sigma(n)$, the subgroup $\KK$ is generated by
$(1,\cdots,1) \in \ZZ^{n+1}$. Thus
$D(\Sigma) \subset (\CC^*)^N$ consists of the elements
$(t,\cdots,t)$, where $t \in \CC^*$. So the toric manifold
associated with $\Sigma(n)$ is the ordinary $n$-dimensional
projective space.
\end{example}

There exists a simple open coverings of $U(\Sigma)$ by affine
algebraic varieties.

\begin{prop}
Let $\sigma$ be a $k$-dimensional cone in $\Sigma$ generated by
$\{v_{i_1},\cdots,v_{i_k}\}.$ Define the open subset $U(\sigma)
\subset \CC^N$ as
$$ U(\sigma) = \{(z_1,\cdots,z_N) \in \CC^N \mid z_j \neq 0
\;\;\textrm{for all}\; j \notin \{i_1,\cdots,i_k\}\}.$$
Then the open sets $U(\sigma)$ have the following properties:
\begin{enumerate}
\item $$U(\Sigma) = \cup_{\sigma \in \Sigma} U(\sigma);$$
\item if $\sigma \prec \sigma'$, then $U(\sigma) \subset U(\sigma')$;
\item for any two cone $\sigma_1,\sigma_2 \in \Sigma$,
one has $U(\sigma_1) \cap U(\sigma_2) = U(\sigma_1 \cap \sigma_2)$;
in particular,
$$ U(\Sigma) = \sum_{\sigma \in \Sigma^{(n)}} U(\sigma).$$
\end{enumerate}
\end{prop}

\begin{prop}\label{homocord}
Let $\sigma$ be an $n$-dimensional cone in $\Sigma^{(n)}$
generated by $\{v_{i_1},\cdots,v_{i_n}\}$, which spans the lattice
$N$. We denote the dual $\ZZ$-basis of the lattice $M$ by
$\{u_{i_1},\cdots,u_{i_n}\}$. i.e.
\begin{equation}
\langle v_{i_k},u_{i_l} \rangle  = \delta_{k,l}
\end{equation}
where $\langle \cdot,\cdot \rangle $ is the canonical pairing
between lattices $N$ and $M$.

Then the affine open subset $U(\sigma)$ is isomorphic to
$\CC^n \times (\CC^*)^{N-n}$, the action of
$D(\Sigma)$ on $U(\sigma)$ is free, and the space of
$D(\Sigma)$-orbits is isomorphic to the affine
space $U_{\sigma} = \CC^n$ whose coordinate functions
$x_1^\sigma,\cdots,x_n^\sigma$ are $n$ Laurent monomials
in $z_1,\cdots,z_N$:
\begin{equation}\label{homocordeq}
\begin{cases}
x_1^\sigma  = z_1^{\langle v_1,u_{i_1} \rangle }\cdots
z_N^{\langle v_N,u_{i_1} \rangle }\\
\qquad \vdots\\
x_n^\sigma  = z_1^{\langle v_1,u_{i_n} \rangle }\cdots
z_N^{\langle v_N,u_{i_n} \rangle }
\end{cases}
\end{equation}

\end{prop}

The last statement yields a general formula for the local affine
coordinates $x_1^\sigma, \cdots,x_n^\sigma$ of a point
$p \in U_{\sigma}$ as functions of its
``homogeneous coordinates'' $z_1,\cdots,z_N$.

\section{Symplectic forms of toric manifolds}
\label{sec:forms}

In the last section, we associated a compact manifold $X_{\Sigma}$
to a fan $\Sigma$. In this section, we review the construction of
symplectic (K\"ahler) manifold associated to a convex polytope
$P$.

Let $M$ be a dual lattice, we consider a convex polytope $P$ in $M_{\RR}$
defined by
\begin{equation}
\{x \in M_{\RR} \mid \langle x,v_j \rangle  \geq \lambda_j
\;\textrm{for}\; j=1,\cdots,N\}
\end{equation}
where $\langle \cdot,\cdot \rangle $ is a dot product of $M_{\RR}
\cong \RR^n$. Namely, $v_j$'s are inward normal vectors to the
codimension 1 faces of the polytope $P$. We associate to it a fan
in the lattice $N$ as follows: With any face $\Gamma$ of $P$, fix
a point $m$ in the (relative) interior of $\Gamma$ and define
$$\sigma_{\Gamma} = \cup_{r \geq 0} r \cdot (P-m).$$
The associated fan is the family $\Sigma(P)$ of dual convex cones
\begin{eqnarray}
\check{\sigma}_{\Gamma} &=&\{ x\in N_\RR \mid \langle  y,x \rangle
\geq 0 \;\;\forall y \in
\sigma_{\Gamma} \} \\
&=&\{ x\in N_\RR \mid \langle  m,x \rangle  \leq \langle p,x
\rangle \;\;\forall p \in P, m \in \Gamma \}
\end{eqnarray}
where $\langle \cdot,\cdot \rangle $ is dual pairing $M_\RR$ and
$N_\RR$. Hence we obtain a compact toric manifold $X_{\Sigma(P)}$
associated to a fan $\Sigma(P)$.

Now we define a symplectic (K\"ahler) form on $X_{\Sigma(P)}$ as
follows. Recall the exact sequence :
$$ 0 \to \KK \stackrel{i}\to \ZZ^N \stackrel{\pi}{\to} \ZZ^n \to 0.$$
It induces another exact sequence :
$$ 0 \to K  \to \RR^N/\ZZ^N \to \RR^n/\ZZ^n \to 0.$$
Denote by $k$ the Lie algebra of the real torus $K$. Then we have
the exact sequence of Lie algebras:
$$ 0 \to k \to \RR^N \stackrel{\pi}{\to} \RR^n \to 0.$$
And we have the dual of above exact sequence:
$$ 0 \to (\RR^n)^* \to (\RR^N)^* \stackrel{i^*}{\to} k^* \to 0.$$

Now, consider $\CC^N$ with symplectic form $\frac{i}{2} \sum dz_k \wedge
d\overline{z}_k$.
The standard action $T^n$ on $\CC^n$ is hamiltonian with moment map
\begin{equation}
\mu(z_1,\cdots,z_N) = \frac{1}{2}(|z_1|^2,\cdots,|z_N|^2).
\end{equation}

For the moment map $\mu_K$ of the $K$ action is then given by
$$
\mu_K=i^* \circ \mu : \CC^N \to k^*.
$$
If we choose a $\ZZ$-basis
of $\KK \subset \ZZ^N$ as
$$ Q_1 = (Q_{11},\cdots,Q_{N1}),\cdots, Q_k= (Q_{1k},\cdots,Q_{Nk})$$
and $\{q^1, \cdots, q^k\}$ be its dual basis of $\KK^*$. Then the
map $i^*$ is given by the matrix $Q^t$ and so we have
\begin{equation}
\mu_K(z_1,\cdots,z_N) = \frac{1}{2} (\sum_{j=1}^N
Q_{j1}|z_j|^2,\cdots,\sum_{j=1}^N Q_{jk}|z_j|^2) \in \RR^k \cong
k^*
\end{equation}
in the coordinates associated to the basis $\{q^1, \cdots, q^k\}$.
We denote again by $\mu_K$ the restriction of $\mu_K$ on
$U(\Sigma) \subset \CC^N$.

\begin{prop}[Audin \cite{audin:ttasm91}, Proposition 6.3.1.]
Then for any $r = (r_1, \cdots, r_{N-n}) \in \mu_K(U(\Sigma))
\subset k^* $, we have a diffeomorphism
\begin{equation}
\mu_K^{-1}(r)/K \cong U(\Sigma)/D(\Sigma) = X_{\Sigma}
\end{equation}
And for each (regular) value of $r \in k^*$, we can associate a
symplectic form $\omega_P$ on the manifold $X_\Sigma$ by
symplectic reduction \cite{MW}.
\end{prop}

To obtain the original polytope $P$ that we started with, we need
to choose $r$ as follows: Consider $\lambda_j$ for $j=1,\cdots,N$
which we used to define our polytope $P$ by the set of
inequalities $\langle x,v_j \rangle  \geq \lambda_j$. Then, for
each $a=1,\cdots,N-n$, let
$$r_a = -\sum_{j=1}^N Q_{ja} \lambda_j.$$
Then we have
$$\mu_K^{-1}(r_1, \cdots, r_{N-n})/K \cong
X_{\Sigma(P)}$$ and for the residual $T^n \cong T^N/K$ action on
$X_{\Sigma(P)}$, and for its moment map $\mu_T$, we have
$$\mu_T(X_{\Sigma(P)}) = P.
$$
In fact, Guillemin \cite{Gu} proved the following explicit closed
formula for the K\"haler form

\begin{theorem}[\bf Guillemin]\label{gullemin}
Let $P$, $X_{\Sigma(P)}$, $\omega_P$ and
$$
\mu_T: X_{\Sigma(P)} \to (\RR^N/k)^*\cong (\RR^n)^*
$$
be the moment map defined as above. Define the functions on
$(\RR^n)^*$
\begin{eqnarray}
\ell_i(x) & = &\langle x, v_i \rangle - \lambda_i \, \mbox{ for }\,
i = 1, \cdots, N \label{eq:elli}\\
\ell_\infty(x) & = &\sum_{i=1}^N\langle x, v_i \rangle =\langle x,
\sum_{i=1}^N v_i \rangle \label{eq:ellinfty}.
\end{eqnarray}
Then we have
\begin{equation}\label{eq:omegaP}
\omega_P =
\sqrt{-1}\partial\overline{\partial}\mu_{T}^*\Big(\sum_{i=1}^N\lambda_i
(\log \ell_i) + \ell_\infty\Big)
\end{equation}
on $\mbox{int}(P)$.
\end{theorem}

\section{Adapted Floer cohomology of the torus fibers}
\label{sec:floer}

Let $(X_{\Sigma(P)},\omega_P)$ be a $2n$-dimensional symplectic
toric manifold with $T^n$-action constructed from the polytope $P
\subset M_\RR$. Each $T^n$ orbit associated to an interior point
in $P$ is a Lagrangian submanifold of $X_{\Sigma(P)}$. Such an
orbit can be obtained as $\mu_T^{-1}(A)$ for $A \in
\mbox{int}(\mu_T(X_{\Sigma(P)}))$ for the moment map $\mu_T$.

We fix one such orbit (non-singular) and denote it by $L$. In this
paper, we will study the Floer cohomology of these Lagrangian tori
and compute this by computing its Bott-Morse theory version
$HF^{BM}(L;J_0)$ as in \cite{cho:hsf03}. One important difference
between the Clifford torus and the general torus fibers is that
the former is {\it monotone} \cite{Oh1}, \cite{Oh4} while the
latters are not. {\it Since the obstruction classes defined in
\cite{FOOO} do not vanish for the Lagrangian submanifold $L$}, it
is not clear whether the standard Floer cohomology $HF(L,\phi(L))$
is defined and invariant under the change of Hamiltonian isotopy,
or whether it is isomorphic to the Bott-Morse version
$HF^{BM}(L;J_0)$ when $L$ is not monotone.

In this section, we will define a restricted version of the Floer
cohomology which exploits some special geometry of Lagrangian
torus fibers in the toric manifolds. We will call this {\it
adapted Floer cohomology} and denote it by $HF^{ad}(L;J_0)$.
Important ingredients for the construction of the adapted Floer
cohomology $HF^{ad}(L;J_0)$ are the following three theorems whose
proof will be postponed to the next two sections.
\smallskip

\noindent {\bf [Maslov index formula]} {\it  For a symplectic
toric manifold $X_{\Sigma(P)}$, let $L$ be a Lagrangian $T^n$
orbit. Then the Maslov index of
 any holomorphic disc with boundary lying on $L$ is twice the
sum of intersection multiplicities of the image of the disc with
the codimension 1 submanifolds $V(v_j)$ for $v_j \in \Sigma^{(1)}$
for all $j=1,\cdots,N$.}
\smallskip

\noindent {\bf [Classification theorem]} {\it Any holomorphic map
$w:(D^2,\partial D^2) \to (X_{\Sigma(P)},L)$ can be lifted to a
holomorphic map
$$\widetilde{w}:(D^2,\partial D^2) \to (\CC^N \setminus Z(\Sigma),\pi^{-1}(L))$$
so that each homogeneous coordinates functions
$z_1(\widetilde{w}), \cdots,z_N(\widetilde{w})$ are given by
Blaschke products with constant factors.
$$i.e. \;\; z_j(\widetilde{w}) = c_j \cdot
\prod_{k=1}^{\mu_j}\frac{z-\alpha_{j,k}}{1-\overline{\alpha}_{j,k}z}$$
for $c_j\in \CC^*$ and non-negative integers $\mu_j$ for each
$j=1,\cdots, N$. In particular, there is no non-constant
holomorphic discs of non-positive Maslov indices.} \smallskip

\noindent {\bf [Regularity theorem]} {\it The discs in the
classification theorem are Fredholm regular, i.e., its
linearization map is surjective.}
\medskip

Assuming these theorems for the moment, we proceed construction of
$HF^{ad}(L;J_0)$ of $(L;J_0)$. We denote the standard integrable
complex structure on $X$ by $J_0$. Let $\phi$ be a Hamiltonian
diffeomorphism such that $\phi(L)$ intersects $L$ transversely. We
consider the set of paths $J^\prime: [0,1] \to \CJ_\omega(X)$ with
$$
J^\prime(0) = J_0, \quad J^\prime(1) = \phi_*J_0
$$
denote it by $j_{(\phi,J_0)}$. Similar theorems obviously hold for
the pair $(\phi(L), \phi_*J_0)$ as for $(L;J_0)$. In particular,
there is no non-constant holomorphic discs of non-positive Maslov
indices for the pair $(\phi(L), \phi_*J_0)$ either.

\begin{remark} The set $j_{(\phi,J_0)}$ was considered and played
an important role in \cite{Oh5} in relation to the formulation of
Floer homology of Hamiltonian diffeomorphisms over the mapping
torus of $\phi$. It appears that considering this set of paths
depending on the triple $(L, J_0; \phi)$ enable us to define the
Floer homology of Lagrangian submanifolds in a more flexible way
when the given pair $(L;J_0)$ has some special structure of the
moduli of $J_0$-holomorphic discs attached to $L$ as in our case.
\end{remark}

Now we restrict to the paths $J^\prime \in j_{(\phi,J_0)}$ for the
study of Floer's equations
\begin{equation}\label{eq:Floer}
\begin{cases}
\frac{\partial u}{\partial \tau} + J^\prime_t \frac{\partial
u}{\partial t} = 0 \\
u(\tau,0) \in L, \, u(\tau,1) \in \phi(L)
\end{cases}
\end{equation}
in the definition of the Floer boundary operator. Now for given
pair $x, y \in L \cap \phi(L)$, we study the moduli space
$$
\CM(x,y;J^\prime)
$$
for the Fredholm index $\mu(x,y) =0, \, 1$ or $2$. The following
proposition is the reason why we restrict $J^\prime$ to the ones
coming from $j_{(\phi,J_0)}$.

\begin{prop}\label{moduli} Assume $X_{\Sigma(P)}$ is Fano.
Let $\phi$ be a Hamiltonian diffeomorphism such that $\phi(L)$
intersects $L$ transversely and let $J^\prime \in j_{(\phi,J_0)}$.
Assume that $x, \, y \in L \cap \phi(L)$ with $\mu(x,y) =0, \, 1$
or $2$. Then the following holds:
\begin{enumerate}
\item When $\mu(x,y) = 0$, $\CM(x,y;J^\prime)/\RR$ is empty.
\item When $\mu(x,y) = 1$, $\CM(x,y;J^\prime)/\RR$ is a compact
manifold of dimension zero
\item When $\mu(x,z) = 2$,
\begin{enumerate}
\item if $x \neq z$, $\CM(x,z;J^\prime)/\RR$ can be compactified into
a compact manifold with boundary of dimension one, whose boundary
consists of the form
\begin{equation}\label{eq:u1u2}
v_1 \sharp v_2
\end{equation}
where $v_1 \in \CM(x,y;J^\prime)$ and $v_2 \in \CM
(y,z;J^\prime)$.
\item if $x =z$, $\CM(x,x;J^\prime)$ can be compactified into a
compact manifold with boundary of dimension one, whose boundary
consists of the types
$$
v_1 \sharp v_2
$$
where $v_i$'s are types either of (\ref{eq:u1u2}) or that for
which one of $v_i$'s is constant and the other is a
$J_0$-holomorphic disc with boundary lying on $L$ or a
$\phi_*(J_0)$-holomorphic disc with boundary lying on $\phi(L)$.
\end{enumerate}
\end{enumerate}
\end{prop}
\begin{proof} First note that non-constant holomorphic discs with boundary on
one of the Lagrangian submanifold, $L$ or $\phi(L)$, have positive
Maslov indices (and so greater than or equal to 2). Once this is
in our disposition, the proof of this fact follows by the
dimension counting arguments from \cite{Oh1}, \cite{Oh4}. We omit
the details of the argument referring to \cite{Oh1}.
\end{proof}
\begin{remark}\label{virtual}
Unlike the case \cite{Oh1} or \cite{Oh4} where we allow to vary
the almost complex structures, since we prefer to keep the usage
of integrable complex structure $J_0$, we also need to prove that
the above singular curves are also regular (or more precisely the
relevant evaluation maps are transverse in forming the fiber
products). This follows from the fact that $L$ is a torus orbit of
the torus action on $X_{\Sigma(P)}$.
\end{remark}

\begin{corollary} Under the hypothesis as in Proposition \ref{moduli},
the Floer cohomology $HF(L, \phi(L);J^\prime)$ is well-defined.
\end{corollary}

We can now compare two Floer cohomology $HF(L,\phi(L);J^\prime)$
with $J^\prime \in j_{(\phi,J_0)}$ and
$HF^*(L,\psi(L));J^{\prime\prime})$ with $J^{\prime\prime} \in
j_{(\psi,J_0)}$ by considering paths
\begin{eqnarray*}
\Phi & = & \{\phi^s\}_{0\leq s \leq 1};\, \phi^0 = \phi, \, \phi^1
= \psi \\
\overline J & = &\{J^s\}_{0\leq s\leq 1}; \, J^0 = J^\prime, \,
J^1 = J^{\prime\prime}, \, J^s \in j_{(\phi^s,J_0)}
\end{eqnarray*}
and the continuity equation
$$
\begin{cases}
\frac{\partial u}{\partial \tau} + J^{\rho(\tau)}_t \frac{\partial
u}{\partial t} = 0 \\
u(\tau,0) \in L, \, u(\tau,1) \in \phi^{\rho(\tau)}(L)
\end{cases}
$$
where $\rho: \RR \to [0,1]$ is a monotonically increasing function
$$
\rho = \begin{cases} 0 \quad \mbox{for }\, \tau \leq -R \\
1 \quad \mbox{for } \, \tau \geq R
\end{cases}
$$
for some sufficiently large $R  >  0$. Again by the same reasoning
using the choice $J^s \in j_{(\phi^s,J_0)}$, we can prove that the
continuity equation defines a chain map
$$
h_{(\Phi,\overline J)}: CF(L,\phi(L); \delta_{J^\prime}) \to
CF(L,\psi(L); \delta_{J^{\prime\prime}})
$$
which is an isomorphism. We refer to \cite{Oh1}, \cite{Oh4} for
the proof in the monotone case, which obviously generalizes in the
current Fano toric case {\it if we use the set-up of the adapted
Floer cohomology}. More specifically we use the special property
of the pair $(L,J_0)$ mentioned in the three theorems in the
beginning of this section. This proves the well-definedness and
the invariance property of $HF(L,\phi(L);J^\prime)$. We denote the
canonical isomorphism class of $HF(L,\phi(L);J^\prime)$ over
$\phi$ and $J^\prime \in j_{(\phi,J_0)}$ by $HF^{ad}(L;J_0)$.

We will compute this group by computing the Bott-Morse version of
the Floer cohomology, which we denote by $HF^{BM}(L;J_0)$. Because
the above structure theorems, this latter Floer cohomology group
is well-defined. The following theorem permits us to do this for
the computation of $HF(L,\phi(L);J^\prime)$.

\begin{theorem}\label{isomorphism}
Assume $X_{\Sigma(P)}$ is Fano and let $L$ and $\phi$, $J^\prime$
as above. Then $HF^{BM}(L;J_0)$ is well-defined and isomorphic to
$HF^{ad}(L;J_0)$. More specifically, $HF^{BM}(L;J_0)$ is
isomorphic to $HF(L,\phi(L);J^\prime)$ for any Hamiltonian
diffeomorphism $\phi$ with $L$ intersection $\phi(L)$ transversely
and a path $J^\prime \in j_{(\phi,J_0)}$.
\end{theorem}
\begin{proof} The well-definedness of $HF^{BM}(L;J_0)$ follows from
the classification theorem which in particular implies that all
holomorphic discs have positive Maslov indices and are regular, by
examining the details of the construction from [section 7, FOOO].

For the second statement, it is enough to consider the case when
$\phi$ is $C^2$-close to the identity. We refer to \cite{Oh4} for
the proof in the monotone case which obviously generalizes to the
semi-positive case, in particular the Fano case $(X,\omega)$.
\end{proof}

\begin{remark}
We would like to emphasize that in our case where the obstructions
do not vanish, the Bott-Morse version of the Floer cohomology
constructed in [FOOO], a priori, may not be defined and may depend
on the choice of the almost complex structure $J$, even if it is
defined. Because of this, we make the dependence on $J_0$ explicit
in the notation $HF^{BM}(L;J_0)$.
\end{remark}

Having Theorem \ref{isomorphism} in mind, we will compute the
Bott-Morse Floer cohomology group $HF^{BM}(L;J_0)$ in the rest of
the paper.

\section{Index formula and the classification of holomorphic discs}
\label{sec:classify}

In this section, we will prove the Maslov index formula and the
classification theorem mentioned in section \ref{sec:floer}.

Before we state the theorem, we recall that for each generator
$v_j \in \Sigma^{(1)}$, there is associated a codimension 1
subvariety $V(v_j)$. For the principle bundle $(U(\Sigma)
\stackrel{\pi}{\to} X_{\Sigma(P)})$, $\pi^{-1}(V(v_j))$ is defined
by the equation $z_j =0$ in $U(\Sigma)$. \medskip

\begin{theorem}[The Maslov index formula]\label{indexformula}
For a symplectic toric manifold $X_{\Sigma(P)}$, let
$L$ be a Lagrangian $T^n$ orbit. Then the Maslov index of
 any holomorphic disc with boundary lying on $L$ is twice the
sum of intersection multiplicities of the image of the disc with
the codimension 1 submanifolds $V(v_j)$ for $v_j \in \Sigma^{(1)}$
for all $j=1,\cdots,N$.
\end{theorem}
\begin{proof} As in \cite{cho:hsf03}, we deform a holomorphic disc
$w:(D^2,\partial D^2) \to (X,L)$ near the intersections with
$V(v_j)$'s. It is easy to see that the intersections are discrete
and there are only finitely many of them because of holomorphicity
of the map $w$. Denote by $p_1 \in D^2$ a point in the preimage of
the intersection. i.e. $p_1 \in w^{-1}(\mbox{image}(w) \cap V(v_j))$
for some $j$. We describe how to deform $w$ as a smooth map near
the point $p_1$ and such deformation will be carried out near every
preimages of intersections.

Note that $w(p_1)$ may lie in the intersection of several
$V(v_j)$'s: Denote them by $V(v_{i_1}),\cdots,V(v_{i_k})$. Then,
We have
\begin{equation}
w(p) \in V(v_{i_1}) \cap \cdots \cap V(v_{i_k})
\end{equation}

The fact that $V(v_{i_1}) \cap \cdots \cap V(v_{i_k}) \neq 0$
implies that $\{v_{i_1},\cdots,v_{i_k}\}$ is not a primitive
collection (See Definition \ref{prim}). Since the fan $\Sigma$ is
complete, we may choose $(n-k)$ generators
$v_{i_{k+1}},\cdots,v_{i_n}$ so that $\langle
v_{i_1},\cdots,v_{i_n} \rangle $ defines a $n$-dimensional cone
$\sigma$ in $\Sigma$.

We may consider the map $w$ near $p_1$ as a map into
the affine open set $\CC^n = Spec(\check{\sigma}\cap M)$
as in \cite{fulton:itv}.
More precisely, the coordinate functions of this affine open set
$\CC^n$ is given as in Proposition \ref{homocord}. Denote by
$d_1,\cdots, d_k$ the intersection multiplicities of the map $w$
with $V(v_{i_1}),\cdots,V(v_{i_k})$. In other words, if we
represent the map $w$ in terms of homogeneous coordinates, then
the homogeneous coordinate functions $z_{i_1},\cdots,z_{i_k}$ will
have order of zero $d_1,\cdots,d_k$ at $p_1$ and other homogeneous
coordinate functions are non-vanishing near $p_1$.

As in Proposition \ref{homocord}, let $\{u_{i_1},\cdots,u_{i_n}\}$
be the basis of $M$ dual to $\{v_{i_1},\cdots,v_{i_n}\}$.
$$\langle u_{i_j},v_{i_k} \rangle  = \delta_{j,k}$$
Then, the affine coordinate function $x^\sigma_1$ is
\begin{eqnarray*}
x^\sigma_1 &=& z_1^{\langle v_1,u_{i_1} \rangle }\cdots
z_N^{\langle v_N,u_{i_1} \rangle }\\
&=& C(z) \cdot z_{i_1}^{\langle v_{i_1},u_{i_1} \rangle } \\
&=& C(z) \cdot z_{i_1}
\end{eqnarray*}
where $C(z)$ is a function
nonvanishing near $p_1$. Therefore, the affine coordinate function
$x^\sigma_1$ has order of zero $d_1$ at $p_1$. Similarly,
$x^\sigma_2,\cdots,x^\sigma_k$ have order of zero $d_2,\cdots,d_k$
at $p_1$. For $j > k$, $x^\sigma_j$ is non-vanishing near zero.
We may further assume that $p_1 = 0 \in D^2$.
Then, the map $w$ near $p_1$ can be written in affine coordinates as
$(a_1z^{d_1} + \CO(z^{d_1+1}),\cdots,a_kz^{d_k} + \CO(z^{d_k+1}),
a_{k+1} + \CO(z),\cdots,a_n + \CO(z)).$

Now we are in the same situation as in \cite{cho:hsf03} Theorem
4.1. From now on, we will only sketch the arguments and refer
readers to \cite{cho:hsf03} for details.

We label by $p_2, \cdots, p_m \in D^2$ all the other points whose
image intersect with $V(v_j)$ for some $j$. We find disjoint open
balls $B_{\epsilon}(p_1) \subset D^2$ centered at $p_i$ with fixed
radius $\epsilon$ for sufficiently small $\epsilon$ for all
$i=1,2,\cdots,m$.

Now we smoothly deform the map $w$ inside the ball
$B_{\epsilon}(p_1) $, so that the deformed map $\widetilde{w}$
satisfies
\begin{equation}
\widetilde{w}|_{\partial B_{\epsilon/2}(p_1)} \subset L
\end{equation}
and as a map into the affine open set $\CC^n$ near $p_1$,
the map $\widetilde{w}$ on $B_{\epsilon/2}(p_1)$ is given by
\begin{equation}\label{map}
\Big(\frac{a_1 z^{d_0}}{|a_0|(\frac{\epsilon}{2})^{d_1}},\cdots,
\frac{a_k
z^{d_k}}{|a_k|(\frac{\epsilon}{2})^{d_k}},\frac{a_{k+1}}{|a_{k+1}|},
\cdots,\frac{a_n}{|a_n|}\Big).
\end{equation}

We perform the same kind of deformations for $p_2,p_3,\cdots,p_m$
inside the ball $B_\epsilon(p_2),\cdots,B_\epsilon(p_m)$ and write
the resulting map as $\widetilde{w}$. Over the punctured disc $$
\Sigma = D^2 \setminus (B_\epsilon(p_1)  \cup \cdots
B_\epsilon(p_m)),$$ the deformed map $\widetilde{w}$ does not
intersect with the hyperplanes, and it intersects with the
Lagrangian torus $L$ along the boundaries of the punctured disc.

Since the Maslov index is a homotopy invariant, we have $\mu(w) =
\mu(\widetilde{w})$. Hence, we may compute the Maslov index of the
map $\widetilde{w}$. Note that the boundary $\partial \Sigma$ is
$\partial D^2 \cup ( \cup_{i} \partial B_{\epsilon/2}(p_i))$.

Since the image of the map $\widetilde{w}$ on the boundaries of
the balls $B_{\epsilon/2}(p_i)'s$ lies on the Lagrangian
submanifold $L$, the map $\widetilde{w}:(\Sigma,\partial \Sigma)
\to (X,L)$ satisfies the Lagrangian boundary condition.
Furthermore, since every intersection with the hyperplane occurs
inside the balls $B_{\epsilon/2}$, $\widetilde{w}|_{\Sigma}$ does
not meet the hyperplanes. Hence, it can be considered as a map
into the cotangent bundle of $L$, (If we take out all such
codimension 1 submanifolds $V(v_j)$'s from $X_{\Sigma(P)}$, there
remains $(\CC^*)^n$ which can be considered as the cotangent
bundle of the torus orbit $L$). Therefore we have
\begin{equation}\label{indexzero}
\mu(\widetilde{w}|_{\Sigma}) = 0.
\end{equation}
On the other hand,
the Maslov index of the map $\widetilde{w}|_{\Sigma}$ is given by the sum
of the
Maslov indices along $\partial \Sigma$ after fixing the trivialization.

Now consider the map $\widetilde{w}:D^2 \to X$ and we fix a
trivialization $\Phi$ of the pull-back bundle $\widetilde{w}^*TX$.
It gives a trivialization $\Phi_{\Sigma}$ of the pull-back bundle
$(\widetilde{w}|_{\Sigma})^*TX$ restricted over $\Sigma$. In this
trivialization, it is easy to see that
$$
\mu(\Phi_{\Sigma},\partial D^2) = \mu(\Phi,\partial D^2) =
\mu(\widetilde{w})=\mu(w).
$$
Since the boundary of the balls $B_{\epsilon/2}$ are
oriented in the opposite way, and from the explicit description
(\ref{map}) of the deformed map on the ball $B_{\epsilon/2}(p_i)$,
we have
$$
\mu(\Phi_{\Sigma},\partial B_{\epsilon/2}(p_i))
= - 2(\textrm{sum of intersection multiplicities in $B_{\epsilon/2}(p_i)$}).
$$
From the equation (\ref{indexzero}), we have
$$\mu(w) - 2(\textrm{sum of intersection multiplicities }) =0.$$
\end{proof}

Now, we use this index formula to classify all holomorphic discs
with boundary lying on $L$. It is much easier if we describe them
in terms of ``homogeneous coordinates'' of toric varieties. Namely
we claim that homogeneous coordinate functions of any holomorphic
discs with boundary on $L$ can be written as just Blaschke
products with constant coefficients.

Before we prove the claim, we first need the following Lemma,
\begin{lemma}
Any holomorphic map $w:D^2 \to X_{\Sigma(P)}$ can be lifted to a
holomorphic map
$$\widetilde{w}:D^2 \to (\CC^N \setminus Z(\Sigma)),$$
so that $$\pi \circ \widetilde{w} = w$$
\end{lemma}
\begin{proof}
The fibration $(\CC^N \setminus Z(\Sigma))  \to  X_{\Sigma(P)}$
is a principal $D(\Sigma)$ bundle.
We pull back this bundle over $D^2$, and fix a holomorphic trivialization
 and take a constant section.
\end{proof}

\begin{theorem}[Classification theorem]\label{classify}
Any holomorphic map $w:(D^2,\partial D^2) \to (X_{\Sigma(P)},L)$ can be
lifted to a holomorphic map
$$\widetilde{w}:(D^2,\partial D^2) \to (\CC^N \setminus Z(\Sigma),\pi^{-1}(L))$$
so that each homogeneous coordinates functions
$z_1(\widetilde{w}), \cdots,z_N(\widetilde{w})$ are given by the
Blaschke products with constant factors.
$$i.e. \;\; z_j(\widetilde{w}) = c_j \cdot
\prod_{k=1}^{\mu_j}\frac{z-\alpha_{j,k}}{1-\overline{\alpha}_{j,k}z}$$
for $c_j\in \CC^*$ and non-negative integers $\mu_j$ for each
$j=1,\cdots, N$.
\end{theorem}

\begin{proof}
By the previous lemma, we have a lift $\widetilde{w}:D^2 \to
(\CC^N \setminus Z(\Sigma))$. Suppose the map $w$ meets  the
submanifold $V(v_1)$ at $w(\alpha)$ for $\alpha \in
\mbox{int}(D^2)$. We multiply factor
$\frac{1-\overline{\alpha}z}{z-\alpha}$ to $z_1(\widetilde{w})$
and denote the modified map by $w_1$. Note that the map $w_1$
still satisfies the boundary condition because
$|\frac{1-\overline{\alpha}z}{z-\alpha}|=1$ for $z \in \partial
D^2$. And the intersection multiplicity of $w_1$ with $V(v_1)$ is
one less than that of $w$.

By repeating the process, we may assume that we obtain a map $w_d$
which does not meet $V(v_1)$. Repeat the process for each $V(v_j)$
for $j=1,\cdots,N$. Hence we obtain a holomorphic map
$\widetilde{w}:D^2 \to (\CC^N \setminus Z(\Sigma))$ which does not
meet any codimension 1 submanifolds $V(v_j)$'s. This map has
Maslov index 0 and is contained in any affine open sets $\CC^n$ of
toric variety. It is easy to see that this map is indeed constant.

Hence, we may deduce that homogeneous coordinates of any
holomorphic disc can be written as Blaschke products.
\end{proof}

\begin{remark}
\begin{enumerate}
\item In the case of $\PP^N$, a similar formula was proved in \cite{cho:hsf03}.
\item The Maslov index of $w$ is $\sum_{j=1}^N \mu_j$ by Theorem \ref{indexformula}.
\end{enumerate}
\end{remark}

\section{Fredholm regularity of discs}\label{sec:regularity}

\begin{theorem}[Regularity theorem]\label{regularity}
The discs in Theorem \ref{classify} are Fredholm regular, i.e.,
its linearization map is surjective.
\end{theorem}

We first recall the exact sequence
$$
0 \to \KK \to \ZZ^N \stackrel{\pi}{\to} \ZZ^n \to 0.
$$
This induces the exact sequence of the complex vector space
$$
0 \to \CC^\KK \to \CC^N \stackrel{\pi}{\to} \CC^n \to 0
$$
via tensoring with $\CC$ where $\CC^\KK$ is the $N-n$ dimensional
subspace of $\CC^N$ spanned by $\KK \subset \ZZ^N$. Note that this
exact sequence is equivariant under the natural actions by the
associated complex tori.

Now we explain implication of the existence of the above
equivariant exact sequence on the study of Fredholm property of
holomorphic map
$$
w: (D^2,\partial D^2) \to (X, L)
$$
where $L \subset X$ is a torus fiber $L = \mu^{-1}(\eta), \, \eta
\in P \subset M_\RR$.

We first need some general discussion on the sheaf of holomorphic
sections of bundle pairs $(E,F)$ where $E$ is a complex vector
bundle over $D^2$ and $F$ a real vector bundle over $\partial D^2$
such that $F \otimes \CC$ an identification with $E|_{\partial
D^2}$. We denote by $(\mathcal E, \mathcal F)$ the sheaf of
holomorphic sections of $E$ with boundary values lying in $F$. We
will be interested in the sheaf cohomology of $(\CE,\CF)$ which we
denote by
$$
H^q(D^2,\partial D^2; E,F) = H^q(E,F).
$$
Here the sheaf cohomology functors are the right derived functors
of the global section functor from the category of sheaves of
$(\CO, \CO_\RR)$-modules on $D^2$ to the category of $\RR$
modules, where $(\CO,\CO_\RR)$ is the sheaf of holomorphic
functions on $D^2$ with real boundary values. Denote by
$\CA^0(E,F)$ the sheaf of $C^\infty$ sections of $E$ with boundary
values in $F$, and denote by $\CA^{(0,1)}(E)$ the sheaf of
$C^\infty$ $E$-valued $(0,1)$-forms. The following is easy to
check (see section 3.4 of \cite{katz:egsm01}).

\begin{lemma}\label{fine} The sequence
$$
0 \to (\CE,\CF) \to \CA^0(E,F) \stackrel{\overline \partial} \to
\CA^{0,1}(E) \to 0
$$
defines a fine resolution of $(\CE,\CF)$.
\end{lemma}

From this, it follows that
\begin{eqnarray*}
H^0(E,F) & \cong & \ker \overline{\partial} \\
H^1(E,F) & \cong & \mbox{coker }\overline{\partial}.
\end{eqnarray*}
Next let $(X,L)$ be a pair of K\"ahler manifold $X$ and a
Lagrangian submanifold $L \subset X$. Consider a holomorphic disc
$w: (D^2,\partial D^2) \to (X,L)$ and denote
$$
E = w ^*TX, \quad F = (\partial w)^*TL.
$$
In terms of the sheaf cohomology group $H^q(D^2,\partial
D^2;E,F)$, the surjectivity of the linearization of the disc $w$
is equivalent to the vanishing result
\begin{equation}\label{eq:vanishing}
H^1(D^2,\partial D^2;E,F) = \{0\}.
\end{equation}

Now we restrict to the case of our main interest as in Theorem
\ref{regularity}. Let $w:(D^2,\partial D^2) \to (X,L)$ be a
holomorphic disc obtained in section \ref{sec:classify} and
$\widetilde w: (D^2,\partial D^2) \to (\CC^N, \pi^{-1}(L))$ be the
lifting obtained in Theorem \ref{classify}. From the expression of
$\widetilde w$ in Theorem \ref{classify}, it follows that
$\widetilde w(\partial D^2)$ is contained in a torus orbit of
$(S^1)^{N}$
$$
\widetilde L = (S^1)^N \cdot (c_1, \cdots, c_N) \subset
\pi^{-1}(L) \subset \CC^N.
$$
We denote by
\begin{eqnarray*}
(E,F) & = &(w^*TX,(\partial w)^*TL) \\
(\widetilde E,\widetilde F) & = & (D^2 \times \CC^N,
(\partial \widetilde w)^*(T\widetilde L))) \\
(E_\KK,F_\KK) & = &((\widetilde w)^*(T Orb_{(\CC_*)^\KK}),
(\partial \widetilde w)^*(T Orb_{K}))
\end{eqnarray*}
and by
$$
(\CE, \CF), \quad (\widetilde \CE,\widetilde \CF), \quad
(\CE_\KK,\CF_\KK)
$$
the corresponding sheaves of holomorphic sections

\begin{lemma}\label{exact}
The natural complex of sheaves
\begin{equation}\label{eq:exact}
0 \to (\CE_\KK,\CF_\KK) \to (\widetilde \CE,\widetilde \CF) \to
(\CE,\CF) \to 0
\end{equation}
is exact.
\end{lemma}
\begin{proof} We need to prove the sequence of stalks
$$
0 \to (\CE_\KK,\CF_\KK)_z \to (\widetilde \CE,\widetilde \CF)_z
\to (\CE,\CF)_z \to 0
$$
is exact at each $z \in D^2$. When $z \in \operatorname{Int}D^2$,
this immediately follows from the $\overline{\partial}$-Poincar\'e
lemma. It remains to prove exactness when $z \in
\partial D^2$. We will give details of the proof of
surjectivity of the last map
\begin{equation}\label{eq:surject}
(\widetilde \CE,\widetilde \CF)_z \to (\CE,\CF)_z
\end{equation}
and leave the rest to the readers.

Let $z_0 \in \partial D^2$. By choosing a sufficiently small
neighborhood $U$ of $z_0$, we can holomorphically identify
$(E_\KK,F_\KK)|_U$ with the trivial bundle $(\CC^{N-n},\RR^{N-n})
\to (U,U\cap
\partial D^2)$. By shrinking $U$ if necessary, we may choose
a holomorphic frame
$$
\{f_1, \cdots, f_{N-n}, f_{N-n+1},\cdots, f_{N}\}
$$
of $(\widetilde E ,\widetilde F)$ so that $f_j = e_j, \, 1 \leq j
\leq N-n$ the standard real {\it constant} basis of $\RR^{N-n}
\subset \CC^{N-n}$ and the projections of $\{[f_{N-n+1}], \cdots,
[f_{N}]\}$ defines a holomorphic frame of $E$.

Now let $\eta$ be a given holomorphic section of $E$ defined in a
neighborhood $z \in V \subset \overline V \subset U$ such that
$$
\eta|_{V \cap \partial D^2} \in F.
$$
We can write
$$
\eta = b_{N-n+1} [f_{N-n+1}] +  \cdots + b_{N} [f_{N}]
$$
where $b_j$'s are holomorphic functions on $V$. Then it is obvious
that
$$
\xi_\eta: = b_{N-n+1} f_{N-n+1}  + \cdots + b_{N} f_{N}
$$
defines a holomorphic section of $\widetilde E$ which projects to
$\eta$. However $\xi_\eta$ may not satisfy the boundary condition
$$
(\xi_\eta)|_{V \cap \partial D^2} \subset \widetilde F
$$
and so we need to correct it by adding a suitable holomorphic
section of $(E_\KK,F_\KK) \cong (\CC^{N-n},\RR^{N-n})$. Since
$(\xi_\eta)|_{\partial D^2} \subset F$, there exists a map $g: V
\cap \partial D^2 \to F = \RR^{N-n}$, $g = (g_1, \cdots, g_{N-n})$
such that
\begin{equation}\label{eq:G}
\xi_\eta|_{V \cap \partial D^2}(z) - \sum_{i=1}^{N-n}g_i(z)e_i \in
\widetilde F
\end{equation}
for all $z \in V \cap \partial D^2$.

Now we solve the following Riemann-Hilbert problem for the map $G:
(D^2,\partial D^2) \to (\CC^{N-n},\RR^{N-n})$, $G=(G_1, \cdots,
G_{N-n})$
\begin{equation}\label{eq:RH}
\begin{cases}
\frac{\partial G}{\partial \overline z}  = 0 \\
G(z) = g(z) \, \quad z \in V \cap \partial D^2
\end{cases}
\end{equation}
It is well-known that this equation can be solved (see \cite{Oh3}
for example) on a neighborhood $V^\prime \subset \overline
V^\prime \subset V$ by multiplying a cut-off function $\rho$ such
that
\begin{eqnarray*}
\rho(z) = \begin{cases} 1 \quad \mbox{for } z\in \overline
V^\prime \\
0 \quad \mbox{for } z \quad \mbox{in a neighborhood of }
\partial {\overline V}.
\end{cases}
\end{eqnarray*}
Now it follows that if we define $\xi$
$$
\xi(z) = \xi_\eta(z) - \sum_{i=1}^{N-n}G_i(z)e_i,
$$
it satisfies
$$
[\xi] = [\xi_\eta] \quad \mbox{and } \xi(z) \in \widetilde F_z,\,
z \in V \cap \partial D^2.
$$
This finishes the proof of surjectivity of (\ref{eq:surject}).
\end{proof}

The exact sequence (\ref{eq:exact}) of the sheaves induces the
long exact sequence of cohomology
\begin{eqnarray}\label{eq:long}
0 &\to & H^0(E_\KK,F_\KK) \to H^0(\widetilde E,\widetilde F) \to
H^0(E,F) \longrightarrow \nonumber\\
&\to & H^1(E_\KK,F_\KK) \to H^1(\widetilde E,\widetilde F) \to
H^1(E,F) \to 0.
\end{eqnarray}
Therefore to prove $H^1(D^2,\partial D^2;E,F)= \{0\}$, it is
enough to prove the following lemma

\begin{lemma} $H^1(\widetilde E, \widetilde F) = \{0\}$.
\end{lemma}
\begin{proof}
From the definition of the bundle pair $(\widetilde E, \widetilde
F) \to (D^2,\partial D^2)$, we have
$$
\widetilde E = D^2 \times \CC^N, \quad \widetilde F = \ell_1
\oplus \cdots \oplus \ell_N.
$$
Here for each $j = 1, \cdots, N$, $\ell_j$ is the line bundle
which is the tangent space of the circle
$$
\theta \mapsto e^{2\pi \mu_j \theta}\cdot c_j \subset \CC
$$
with $\mu_j \geq 0$ is an integer given in Theorem \ref{classify}.
Now the lemma immediately follows from the study of the
one-dimensional Riemann-Hilbert problem with this Lagrangian loop
(see e.g.,  \cite{Oh3}  for this kind of analysis).
\end{proof}

This finishes the proof of the vanishing result
$$
H^1(E,F) = \{0\}
$$
and so the discs $w$ obtained in Theorem \ref{classify} and so all
the discs in $X$ with boundary lying on $L$ are Fredholm-regular.

\section{Holomorphic discs of Maslov index two}\label{sec:index-two}

We first recall the definition of the Bott-Morse Floer coboundary
operator from \cite{FOOO}: For $[P,f]\in C^*(L,\QQ)$ and non-zero
$\beta \in \pi_2(M,L)$,
\begin{equation}\label{del}
\begin{cases}
 \D_{\beta}([P,f])=(\CM_2(\beta) \,_{ev_1}\times_f P, ev_0) \\
 \delta_0 ([P,f]) = (-1)^n [\partial P, f]
\end{cases}
\end{equation}
And the boundary operator is defined as
\begin{equation}\label{delall}
\delta([P,f])=\sum_{\beta \in \pi_2(M,L)} \D_\beta([P,f]) \otimes
T^{\omega(\beta)}q^{\frac{\mu(\beta)}{2}}
\end{equation}
And we extend it linearly over the universal Novikov ring
$\Lambda_{nov}$. The following boundary property follows from the
proof of [Theorem 6.24, FOOO] in which is considered the case
where all the obstructions vanish, after combined with some
additional cancellation arguments used in [addenda, O1], [Theorem
2.28, Cho] to deal with the case where the obstruction does not
vanish but is a multiple of the fundamental cycle. We omit the
proof referring to that of [Theorem 2.28, Cho].

\begin{theorem} Assume that $X_{\Sigma(P)}$ is Fano and $L$ is as before.
Then
$$
\delta\circ \delta = 0.
$$
\end{theorem}

Since the standard complex structure $J_0$ in these toric
manifolds are regular as proved in the last section, we may
proceed to compute the actual Floer boundary map with respect to
$J_0$. The relevant calculations in our cases will be reduced to
the study of discs of Maslov index two  as in \cite{cho:hsf03}
because of the following proposition.

\begin{prop}\label{zero4} Let $\delta_k$ to be the formal sum of
$\delta_\beta$ with $\mu(\beta) = k$. Then we have $\delta_k
\equiv 0$ for $k \geq 4$.
\end{prop}
\begin{proof} We can proceed as in the case of Clifford torus. Consider
the homotopy class $\beta \in \pi_2(X,L)$ with the Maslov index
$\mu(\beta)=  4$. The fiber product in the Floer coboundary
operator
 $\CM(\beta)_{ev_1} \times_f P$ has expected dimension
 $dim(P) +3$. If its dimension is less than expected dimension,
the boundary operator is considered as zero since we consider them
in terms of currents (See \cite{FOOO} for details).
 But it is not hard to see that the dimension of
$\CM(\beta)_{ev_1} \times_f P$ is less or equal to $dim(P) +2$ :
Consider the case that $P$ is a point cycle $\langle pt \rangle $
in $X$. The fiber product $\CM(\beta)_{ev_1} \times_f P$ is
nothing but the image of the trajectories of the boundary of
holomorphic discs in $\CM(\beta)$ which meets the point $\langle
pt \rangle $. Consider the lifts of these holomorphic discs in
$U(\Sigma)$. Then from the expression of lifted discs in Theorem
\ref{classify}, the image of the boundary of the lifted discs has
dimension always less or equal to two which is the dimension of
$S^1\times S^1$. Hence after taking a quotient by $D(\Sigma)$,
dimension of its image is still less than two. This proves the
proposition for the case $P$ is a point cycle $\langle pt \rangle
$, and other cases can be done similarly.
\end{proof}

Therefore, we will be mainly interested in the holomorphic discs
of Maslov index two for the computation of $HF^{BM}(L;J_0)$. From
the classification theorem, it is easy to see that there exists
$N$ number of  holomorphic discs of Maslov index 2 (up to an
automorphism of a disc) meeting a point in $L$. We denote the
homotopy class of such discs by $\beta_j \in \pi_2(X,L)$ for
$j=1,\cdots,N$:
\begin{definition}
For the homogeneous coordinates $z_1,\cdots,z_N$, we denote by
$D(v_j)$ the holomorphic disc of class $\beta_j \in \pi_2(X,L)$
associated to the lifted disc
\begin{equation}\label{dj}
\begin{cases}
z_k =c_k  \;\;\textrm{for} \;\; k \neq j \\
z_j = c_j \cdot z
\end{cases}
\end{equation}
for $z \in D^2$, where $(c_1,\cdots,c_N) \in (\CC^*)^N$ are  chosen
to satisfy the boundary condition.
\end{definition}

Now we want to express each such disc in terms of the coordinates
of the torus $(\CC^*)^n \subset X_{\Sigma(P)}$ to compute the
boundary operator.
Recall that in toric varieties, the torus $(\CC^*)^n$ corresponds
to $0$-cone in $N$ or the dual cone $M_\RR$.
$$
(\CC^*)^n \cong Spec \; \CC[x_1,x_1^{-1},x_2,x_2^{-1},\cdots,x_n,x_n^{-1}].
$$
Its coordinate can also be obtained by applying Proposition
\ref{homocord} for the cone $\sigma$ which is generated by the
standard basis vectors $\langle e_1,\cdots,e_n \rangle $ (Such
cone may not exist in the fan $\Sigma$, but the coordinate
expression of $(\CC^*)^n$ obtained this way is still true).

Hence we use Proposition \ref{homocord} to find the relation with
the $(\CC^*)^n$ coordinates and the homogeneous coordinates. If we
choose the generators of the cone $(v_{i_j})$ in Proposition
\ref{homocord} to be $\langle e_1,\cdots, e_n \rangle $, its dual
basis becomes
$$ u_{i_j} = e_j^*.$$
 From the equation (\ref{homocordeq}), we have
\begin{equation}
\begin{cases}
x_1^\sigma = z_1^{\langle v_1,e_1^* \rangle }\cdots z_N^{\langle v_N,e_1^* \rangle } \\
\qquad \vdots \\
x_n^\sigma  =  z_1^{\langle v_1,e_n^* \rangle }\cdots z_N^{\langle
v_N,e_n^* \rangle }
\end{cases}
\end{equation}
Hence for the holomorphic disc $D(v_j)$, by substituting
(\ref{dj}) into the above equations, we get the following :
\begin{equation}
\begin{cases}
x_1^\sigma  = c_1' \cdot z^{\langle v_j,e_1^* \rangle }
=  c_1' \cdot z^{v_j^1}\\
\qquad \vdots \\
x_n^\sigma  =  c_n' \cdot z^{\langle v_j,e_n^* \rangle } = c_n'
\cdot z^{v_j^n}
\end{cases}
\end{equation}
where $v_j = (v_j^1, \cdots, v_j^n)$.
\begin{prop}\label{disccord}
For $i=1,\cdots,N$, the holomorphic disc $D(v_j)$ given by
(\ref{dj}) can be written in terms of coordinates of the torus
$(\CC^*)^n$ as
\begin{equation}\label{zui}
(C_1 z^{v_j^1},C_2 z^{v_j^2},\cdots,C_n z^{v_j^n})
\end{equation} where
constants $C_i \in \CC$ are chosen to satisfy the given Lagrangian
boundary condition.
\end{prop}

\begin{example}
For the Clifford torus case, the holomorphic discs of index two
are
$$ [z:c_1:\cdots:c_n], \cdots, [1:c_1:\cdots:c_n z],$$
 which in the standard open set $U_0$ are
$$ (c_1\frac{1}{z},\cdots,c_n\frac{1}{z}),(c_1 z,\cdots,c_n),
\cdots,(c_1,\cdots,c_n z)$$ Now, the image of the moment map of
$\PP^n$ is the standard $n$ simplex, which can be written as
follows:

For $v_1 =e_1, v_n=e_n, v_{n+1}=(-1,-1,\cdots,-1) \in \RR^n$,
\begin{equation}
\begin{cases}
\langle x,v_i \rangle  \;\geq \;0 \;\; for \;\; i \leq n \\
\langle x,v_{n+1} \rangle  \;\geq \;-1
\end{cases}
\end{equation}

Now one can see the theorem is true in this case.
\begin{equation}
\begin{cases}
v_{n+1} \Longrightarrow (c_1\frac{1}{z},\cdots,c_n\frac{1}{z}) \\
v_j \Longrightarrow  (c_1,\cdots,c_j z,\cdots,c_n)
\end{cases}
\end{equation}

\end{example}

We have the classification theorem, Theorem \ref{classify} in
terms of the homogeneous coordinates, but it is also convenient to
look at them in the open sets $\CC^n$ corresponding to
$n$-dimensional cones in $\Sigma$. But one should note that not
all discs are contained in these affine open sets. More precisely,
if the holomorphic disc intersects with
$V(v_{i_1}),\cdots,V(v_{i_j})$ (possibly at different points), and
if $\{v_{i_1},\cdots,v_{i_j}\}$ is a {\em primitive collection},
then such disc can not be contained in the affine open sets. But
as the primitive collections have two or more elements, the discs
of Maslov index two  which intersect only one of the submanifolds
$V(v_j)$'s are always contained in the affine open sets.

\begin{prop}
For the affine open set $\CC^n$ corresponding to $n$-dimensional
cone $\sigma = \langle  v_{i_1},\cdots,v_{i_n} \rangle $ in
$\Sigma$, the holomorphic discs with Maslov index 2 contained in
this open set $\CC^n \subset X$ are just
$D(v_{i_1}),\cdots,D(v_{i_n})$ up to an automorphism of a disc.
\end{prop}
\begin{proof}
For such an open set $\CC^n \subset X$, the Lagrangian torus fiber
$L$ is defined by $|z_i| =c_i$ for $i=1,\cdots,n$ for some $c_i
\in \RR$. And the holomorphic discs which are mapped into this
open set $\CC^n$ are indeed easy to classify. More precisely, the
$i$-th coordinate of such maps are just given by the Blaschke
products times the constant $c_i$. Hence, holomorphic discs of
Maslov index 2 are (up to automorphism of disc) can be written in
terms of coordinates of $\CC^n$ as
\begin{eqnarray*}
&(c_1  z , c_2 , \cdots , c_n)& \\
&(c_1 , c_2  z , \cdots , c_n)& \\
&\vdots& \\
&(c_1 , c_2 , \cdots , c_n  z)&
\end{eqnarray*}

As the coordinate of $\CC^n$ is determined by the dual cone
$\check{\sigma}$ of the cone $\sigma = \langle
v_{i_1},\cdots,v_{i_n} \rangle $. The primitive generators of
$\check{\sigma}$ are given by the dual $\ZZ$-basis $\langle
u_1,\cdots,u_n \rangle $ in $M$ since $X$ is smooth.

Let $z_1,\cdots,z_n$ be the coordinates of the torus $(\CC^*)^n
\subset X$ given by $M_\RR$. From \cite{fulton:itv}, the affine
coordinates $x^\sigma_1,\cdots,x^\sigma_n$ are given by the
primitive generators as follows: For $u_i :=(u_{i1},\cdots,u_{in})
\in M$,
\begin{equation}
\begin{cases}
x^\sigma_1 & = z_1^{u_{11}}z_2^{u_{12}}\cdots z_n^{u_{1n}} \\
&\vdots \\
x^\sigma_n &=  z_1^{u_{n1}}z_2^{u_{n2}}\cdots z_n^{u_{nn}}
\end{cases}
\end{equation}

Then, the torus coordinates $z_k$ can be recovered from the affine
coordinates $x^\sigma_1,\cdots,x^\sigma_n$: Take
\begin{eqnarray*}
(x^\sigma_1)^{v_{i_1}^k} \cdot (x^\sigma_2)^{v_{i_2}^k} \cdots
(x^\sigma_n)^{v_{i_n}^k} &= &z_1^{(u_{11}v_{i_1}^k + \cdots +
u_{n1}v_{i_n}^k)} \cdots z_n^{(u_{1n}v_{i_1}^k + \cdots +
u_{nn}v_{i_n}^k)} \\
& = & z_1^{(V^t\cdot U)_{1k}} \cdots z_n^{(V^t\cdot U)_{nk}} = z_k
\end{eqnarray*}
where $U$, $V$ are $(n\times n)$ matrices whose $j$-th rows are
given by the vectors $v_{i_j}$, $u_j$ respectively. The last
equality follows from the duality between $v_{i_j}$ and $u_j$.

 Hence the holomorphic disc in $\CC^N$ given by
$$(c_1 ,\cdots,c_j z,\cdots , c_N) $$
can be rewritten in the coordinates of the torus $(\CC^*)^n$ as
\begin{equation}
\begin{cases}
z_1 & = (x_1^\sigma)^{v_{i_1}^1}\cdots (x_n^\sigma)^{v_{i_n}^1}
= C_1 \cdot z^{v_{i_j}^1} \\
&\vdots \\
z_n & = (x_1^\sigma)^{v_{i_1}^n}\cdots (x_n^\sigma)^{v_{i_n}^n}
= C_n \cdot z^{v_{i_j}^n} \\
\end{cases}
\end{equation}
for $(C_1,\cdots,C_n) \in (\CC^*)^n$. This is nothing but the
expression of the disc $D(v_{i_j})$ in Proposition \ref{disccord}.
This proves the proposition.
\end{proof}

\section{The areas of holomorphic discs}
\label{sec:area}

In this section we compute the symplectic areas of the holomorphic
discs. For each such holomorphic disc $D(v_j)$, there exists
$S^1$-action on its image from the torus action on the toric
variety. From the coordinate expression of holomorphic discs in
Theorem \ref{disccord}, this $S^1$ can be easily seen as a
subgroup of $T = (S^1)^n$ via the {\it monomorphism}
\begin{equation}\label{s1action}
S^1 \to T : e^{i\theta} \mapsto (e^{iv_j^1\theta},
\cdots,e^{iv_j^n\theta})
\end{equation}
for each given $j= 1, \cdots, N$. We will fix one such $j$ in the
rest of this section.

 In the level of Lie
algebra, the $S^1 \subset T$ is generated by the element
\begin{equation}\label{eq:xi}
\xi = v_j^1e_1 + v_j^2 e_2 +  \cdots + v_j^n e_n \in Lie(T^n)
\cong \RR^n
\end{equation}

From now on, we denote by $\mu_T$ for the moment map of the whole
torus $(T \cong (S^1)^n)$ action. The image $\mu_T(D(v_j))$ of
holomorphic discs $D(v_j)$ under the moment map $\mu_T$ can be
easily seen to be 1-dimensional because it is invariant under the
$S^1$ action generated by $\xi$, and it meets with the boundary of
the moment polytope because when the disc meets the submanifold $V(v_j)$.
The intersection point is a fixed point of
the $S^1$ action we described above. Indeed, $\mu_T(D_j)$ meets
the hyperplane defined by
$$\langle x, v_j \rangle\;= \; \lambda_j,$$ since the preimage
under the moment map $\mu_T$ of this hyperplane has the stabilizer
$v_j$. Also recall that the image of the Lagrangian torus fiber
under $\mu_T$ is a point, which we denote by
$$A=(a_1,a_2,\cdots,a_n) \in (\RR^n)^*.$$

Let $(r,\theta)$ be the standard polar coordinate of
$D^2(1)\subset \CC$ and consider the map
$$
(r,\theta) \mapsto \mu_T(w(r,\theta))
$$
where
$$
w=D(v_j) :(D^2,\partial D^2) \to (X,L)
$$
provided in
Proposition \ref{disccord}. Since the disc is invariant under the
$S^1$-action (\ref{s1action}), the map is independent of $\theta$.
We write the corresponding curve by
$$
\alpha: [0,1] \to (\RR^n)^*=(Lie(T^n))^*; \quad \alpha(r):=
\mu_T(w(r,\cdot))
$$
We are now ready to prove the following area formula of the disc
$D(v_j)$, which will play a crucial role later when we relate our
computation of the Floer cohomology to Hori-Vafa's Landau-Ginzburg
$B$-model calculation.

\begin{theorem}\label{area}
The area of the holomorphic disc $D(v_j)$ in Proposition
\ref{disccord} is
$$2\pi(\langle A, v_j\rangle  - \lambda_j).$$
\end{theorem}
\begin{proof} Let $\eta \in Lie(T^n)$ be any element and $\eta_X$ be the
vector field on $X$ generated by $\eta$. By definition of the
moment map $\mu_T$, we have the following defining formula of the
moment map
$$
d\langle \mu_T,  \eta \rangle = \eta_X \rfloor \omega_P
$$
in general \cite{MW}. We apply this identity to $\eta = \xi$
defined in (\ref{eq:xi}) to have
\begin{equation}\label{eq:dmuT}
d\langle \mu_T, \xi \rangle = \xi_X \rfloor \omega_P.
\end{equation}
Therefore we derive
\begin{eqnarray}
\frac{d}{dr}\langle \alpha(r), \xi \rangle & = & \langle
d\mu_T\Big(\frac{\partial w}{\partial r}\Big), \xi \rangle \nonumber \\
& = & d \langle \mu_T, \xi \rangle \Big(\frac{\partial w}{\partial
r}\Big) \nonumber \\
& = & \xi_X \rfloor \omega_P \Big(\frac{\partial w}{\partial
r}\Big) \label{eq:xiT}
\end{eqnarray}
where we regard $\mu_T$ both as the map from $X$ to $(Lie(T^n))^*$
and as a $(Lie(T^n))^*$-valued function.

And it follows from the coordinate formula (\ref{zui}) that
$$
\xi_X(w(r,\theta)) = \frac{\partial w}{\partial \theta}(r,\theta).
$$
By substituting this into (\ref{eq:xiT}), we have derived
\begin{equation}\label{eq:ddralpha}
\frac{d}{dr}\langle \alpha(r), \xi \rangle =
\omega_P\Big(\frac{\partial w}{\partial \theta}, \frac{\partial
w}{\partial r}\Big).
\end{equation}
From this, we derive
\begin{eqnarray*}
Area(D(v_j)) & = & \int_{D^2} w^*\omega_P =  \int_0^1
\int_0^{2\pi} \omega_P\Big(\frac{\partial w}{\partial r},
\frac{\partial w}{\partial \theta}\Big) d\theta\, dr \\
& = & - 2 \pi \int_0^1\frac{d}{dr}\langle \alpha(r), \xi \rangle\, dr\\
& = &2\pi(\langle \alpha(0), \xi \rangle - \langle \alpha(1), \xi
\rangle).
\end{eqnarray*}
The value of $\alpha(1) \equiv \mu_T(w(1,\theta))$ is the base of
the Lagrangian torus fiber $L$ which is nothing but $\langle A,
\xi \rangle $ and $\alpha(1)$ is in the hyperplane determined by
$$\langle x, \xi\rangle = \lambda_j.$$
Therefore we have proved that the area of the disc  is
$2\pi(\langle A, \xi\rangle - \lambda_j)$. Finally noting that
$\xi = v_j$ in (\ref{eq:xi}), we have finished the proof.
\end{proof}

\section{Standard spin structure}
We recall the notion of the {\em standard spin structure}
introduced in \cite{cho:hsf03} for the case of the Clifford torus
in $\PP^n$. A spin structure of $L$ is equivalent to the homotopy
class of a trivialization of the tangent bundle of $L$ over the
two skeleton of $L$. We also recall that a {\it framing} of the
manifold $L$ is defined to be the homotopy class of a
trivialization of the tangent bundle $TL$. Therefore each framing
canonically fixes a spin structure of $L$.

\begin{prop}\label{prop:standard}
The framings of $L$ induced by the embeddings $L \hookrightarrow
U_\sigma \cong \CC^n$ are all the same over the choice of cones
$\sigma$. We call the corresponding spin structure of $L$ the {\em
standard spin structure} of $L \subset X$.
\end{prop}
\begin{proof}
Let $S^1:= e^{i\theta}$  be the unit circle embedded in $\CC$. The
tangent bundle of $S^1$ has a natural trivialization given by $S^1
\times \RR \cdot \PD{\theta}$. Similarly there is a natural
trivialization of the tangent bundle of $(S^1)^n \subset \CC^n$.
The torus fiber $L$ sits inside the intersection of the affine
open sets. So, each affine open set induces a trivialization of
tangent bundle of $L$. One can check that the trivializations of
$TL$ obtained for each affine open set have the same homotopy
class because the transition matrices are constant matrices:
Recall that the transition functions between these affine open
sets are given by monomial relations. For two $n$-dimensional
cones $\sigma, \tau$, Let
$$z_i^\sigma = (z_1^\tau)^{a_{i1}}\cdots
(z_n^\tau)^{a_{in}}.
$$
Then,
$$
\frac{\partial}{\partial \theta_{z_j^\tau}}
= a_{1j} \frac{\partial}{\partial \theta_{z_1^\sigma}}
+ \cdots + a_{nj} \frac{\partial}{\partial \theta_{z_n^\sigma}}.
$$
Hence, the transition matrices for the induced trivializations of
the $TL$ are constant matrices. By permuting the affine
coordinates to make $\det(a_{ij}) > 0$, if necessary, this implies
that the trivializations induced by each affine sets are in the
same homotopy class. This is what we mean by the standard spin
structure of $T^n$.
\end{proof}

Recall from [FOOO] that to fix an orientation of the moduli spaces
of holomorphic discs we need to fix a spin structure of $L$ and an
identification of the tangent space at a point of $L$ with
$\RR^n$. Different trivializations in the same homotopy class can
only reverse signs of the moduli spaces of the holomorphic discs
simultaneously and hence give the same Floer cohomology group.

Also note that there exists $2^n = |H^1(L;\ZZ/2)|$ different spin
structures for the torus $L$. Other spin structures besides the
standard one can be naturally considered in the setting of the
Floer cohomology twisted by the flat line bundles on $L$:
Calculations of the Floer cohomology with different spin
structures can be substituted by the Floer cohomology twisted by
the flat line bundles on $L$ with holonomy $e^{\pi i}$ along
appropriate generators of $\pi_1(L)$.

Our computations in the rest of the paper will be based on the
standard spin structure. We refer readers to \cite{cho:hsf03} for
more detailed discussions on the orientation and computations for
different spin structures.

\section{Computation of the Bott-Morse Floer cohomology}
\label{sec:computation}

Now, we are ready to compute the Bott-Morse Floer cohomology of
any Lagrangian torus fiber $L$ in symplectic toric manifold
$X_{\Sigma(P)}$. We will assume in this section that
$X_{\Sigma(P)}$ is Fano.

The Bott-Morse Floer cohomology  defined in section
\ref{sec:index-two} satisfies
$$
\delta \circ \delta =0
$$
for our torus fiber.  Note that we do not need to deform the
boundary operator of the Floer complex by introducing obstruction
cycles since all non-constant holomorphic discs have positive
Maslov indices in our case.

We fix the standard spin structure of $L$, which fixes the
orientation of the moduli space of holomorphic discs. The
orientation of the boundary (\ref{del}) not only depends on the
orientation of the moduli space $\CM_2(\beta)$, but also the fiber
product orientation. It was studied in great detail in
\cite{cho:hsf03}, \cite{FOOO}, and so we restrict our discussion
about orientation to a minimum.

Recall that the Floer cochain complex in \cite{FOOO} is
constructed using currents. From now on, the cycles we write
actually represents their Poincar\'e duals, and we will not
distinguish homology $H_*(L,\QQ)$ and cohomology $H^*(L,\QQ)$ in
our presentation.

The filtration on the boundary operator $\delta$ with energy
induces a spectral sequence $E^{*,*}_r$ which converges to the
Floer cohomology $HF^{BM}(L;J_0)$. Recall from \cite{FOOO} that
$$
E_2^{p,q} \cong (H^*(L,\QQ) \otimes e^q)^p
$$
where $(\;)^p$ means the total degree
$p$. To compute the Floer cohomology, we work with this spectral
sequence and the main step is to compute the boundary $\delta_2$
of the cohomology generators. Here $\delta_2$ is the boundary
operator given by considering only Maslov index 2 discs.

We first compute the boundary for a point class $\langle pt \rangle $.
We denote the generators of $H_*(L,\QQ)$ by $L_1,\cdots,L_n$. More
precisely, by $L_j$ we denote a cycle given by the image of the
map $$S^1 \to (\CC^*)^n : e^{i\theta} \mapsto
(c_1,\cdots,c_je^{i\theta}, \cdots,c_n).$$

In view of Proposition \ref{disccord} and considering an
orientation as in \cite{cho:hsf03}, we have
\begin{equation}
\delta_{\beta_j}\langle pt \rangle  = (-1)^n (v_j^1L_1 + \cdots +
v_j^n L_n)
\end{equation}
where $\beta_j = [D(v_j)] \in \pi_2(X,L)$. Hence,
\begin{eqnarray*}
\delta_2(\langle pt \rangle )&=&\sum_{j=1}^N (-1)^n
T^{Area(\beta_j)}\cdot q\cdot
(v_j^1 L_1 + \cdots + v_j^nL_n)\\
&=&\sum_{j=1}^N (-1)^n T^{2\pi(\langle v_j,A \rangle
-\lambda_j)}\cdot q \cdot (v_j^1L_1 + \cdots + v_j^nL_n)
\end{eqnarray*}
We can also compute the Floer cohomology with flat line bundle
$\CL$ on it, which we denote by
$$
HF^{BM}((L,\CL);J_0).
$$
If we denote by
$$
h_\alpha = e^{i \nu_\alpha}
$$
the holonomy of the line bundle $\CL$ along the cycle $L_\alpha$
for $\alpha = 1, \cdots, n$, Proposition \ref{classify} implies
that the holonomy along the boundary of the disc $D(v_j)$, becomes
\begin{equation}\label{hvj}
h_1^{v_j^1}\cdots h_n^{v_j^n} = e^{i \langle \nu, v_j \rangle} :=
h^{v_j}
\end{equation}
where the vector $\nu = \nu_\CL$ is defined by
\begin{equation}\label{eq:nu}
\nu = (\nu_1, \cdots, \nu_n)
\end{equation}
which we call the {\it holonomy vector} of $\CL$.

In this case, the boundary operator of the Floer cochain complex
is defined as follows \cite{Fuk2}:
\begin{equation}
\begin{cases}
\delta_{\beta}([P,f])=(\CM_2(\beta) \,_{ev_1}\times_f P,
ev_0)\cdot
 (hol_{\partial \beta} \CL) \otimes q \,\,\,\,\textrm{for }\, \beta \neq 0
 \\
\delta_0([P,f])=(-1)^n[\partial P,f]
\end{cases}
\end{equation}
Therefore, we have
\begin{equation}
\delta_2(\langle pt \rangle )=\sum_{j}(-1)^n  h^{v_j}
T^{2\pi(\langle v_j,A \rangle -\lambda_j)}\cdot q \cdot (v_j^1L_1
+ \cdots + v_j^n L_n)
\end{equation}
By identifying $H_1(L:\QQ)$ with $\QQ^n$ via $L_i \mapsto e_i$, we
may write the condition to have $\delta_2(\langle pt \rangle )=0$
as
\begin{equation}\label{holboundary}
\sum_{j}(-1)^n  h^{v_j} T^{2\pi(\langle v_j,A \rangle
-\lambda_j)}\cdot v_j =0
\end{equation}
It is not hard to see that if $\delta_2 (\langle pt \rangle ) =0$,
we would have $\delta_2(P) =0$ in $H^*(L,\QQ)$ for any cycle $P
\in H_*(L,\QQ)$ (see [Cho] for the relevant computations).
Therefore, in this case, the Floer cohomology
$HF^{BM}((L,\CL);J_0)$ is isomorphic to the singular cohomology of
$L$. In particular, it is non-vanishing. The following proposition
implies that, one only needs to consider $\delta_2 (\langle pt
\rangle )$ for the computation of Floer cohomology.

\begin{theorem}\label{criterion}
If $\delta_2 (\langle pt \rangle )=0$, then Bott-Morse Floer
cohomology is isomorphic to the singular cohomology of $L$ as a
$\Lambda_{nov}$-module, i.e.,
$$HF^{BM}((L,\CL);J_0) \cong H^*(L;\Lambda_{nov}^\CC)$$
where $\Lambda_{nov}^\CC$ is the Novikov ring twisted by the line
bundle $\CL$ in an obvious way.

If $\delta_2 (\langle pt \rangle ) \neq 0$, then the Floer
cohomology $HF^{BM}((L,\CL);J_0)$ vanishes.
\end{theorem}
\begin{proof}
 It remains to prove the second statement. Suppose $\delta_2 (\langle pt
\rangle )  \neq 0$, and consider the lowest energy terms of
$\delta_2\langle pt \rangle $ which gives rise to a non-zero term:
Suppose the terms with this energy are given by
$\delta_{\beta_{i_1}},\cdots,\delta_{\beta_{i_\ell}}$. Denote by
$\widetilde{\delta}_2$ the sum
$$\widetilde{\delta}_2 := \delta_{\beta_{i_1}}+ \cdots+\delta_{\beta_{i_\ell}}$$
By the assumption $\delta_2 \langle pt \rangle \neq 0$, we have
$\widetilde{\delta}_2 \neq 0$. It follows from the construction of
the spectral sequence in \cite{FOOO} that this becomes the
boundary operator of the spectral sequence of a certain step, say
$r$. From our choice of $\beta_{i_*}$, lower energy terms give
rise to zero boundary maps in the spectral sequence. Therefore we
have,
$$E_r^{p,q} \cong E_2^{p,q} \cong (H^*(L,\QQ) \otimes e^q)^p$$
We will show that $$E_{r+1} \cong 0.$$ For this we will compute
$\widetilde{\delta}_2$ for the cohomology generators of
$H^*(L,\QQ)$.

In $H^*(L,\QQ)$, we may write, omitting the common factor of
formal parameter $T^{Area}q$,
\begin{equation}
\widetilde{\delta}_2 \langle pt \rangle  = c_1[L_1] + c_2[L_2] +
\cdots + c_n[L_n].
\end{equation}
At least one of $c_i$ is non-zero from our assumption. It is not
hard to see that
$$\widetilde{\delta}_2\langle L_i \rangle  =
\sum_{j=1}^n c_j \langle L_j \times L_i \rangle $$
where $L_i \times L_i$ is 0-cycle. Or more generally,
$$\widetilde{\delta}_2(L_{i_i}\times L_{i_2} \times \cdots \times L_{i_k})=
\sum_{j=1}^n c_j\langle L_j \times (L_{i_i}\times L_{i_2}
\times \cdots \times L_{i_k}) \rangle $$ where the latter is a
0-cycle if $j \in \{i_1,i_2,\cdots,i_k\}$ (See \cite{cho:hsf03}
for the case of Clifford torus in $\PP^n$).

From now on, for index sets, say $J$ with $j=|J|$ elements,
we denote its elements as
$J=\{j_1,\cdots,j_j\}$ with $j_1 < j_2 < \cdots < j_j$.
And we denote $J_{\widehat{s}} = J \setminus \{j_s\}$.

Now we denote an arbitrary element of $k$ dimensional cycles as
$$\sum_{I,|I|=k} A_I L_I$$
for $I \subset \{1,2,\cdots,n\}$ and $A_I \in \QQ$. The boundary
of this element is
\begin{eqnarray*}
\widetilde{\delta}_2(\sum_{I,|I|=k} A_I L_I) &=& \sum_{I} A_I
(\widetilde{\delta}_2L_I) \\
&=& \sum_{I} A_I (c_1L_1 +  \cdots + c_nL_n)\times L_I \\
&=& \sum_{J,|J|=k+1}\sum_{s=1}^{k+1} A_{J_{\widehat{s}}}(-1)^{s-1}
c_{j_s} L_J
\end{eqnarray*}

Hence, the element $\sum_{I,|I|=k} A_I L_I$
is in the kernel of $\widetilde{\delta}_2$ if for any set
$J\subset \{1,2,\cdots,n\}$ with
$|J|=k+1$, the following equation holds:
\begin{equation}\label{kernel}
\sum_{s=1}^{k+1} A_{J_{\widehat{s}}}(-1)^{s-1} c_{j_s}=0.
\end{equation}

Set
\begin{equation}
\begin{cases}
S:= \{ i\in \{1,2,\cdots,n\} \mid c_i =0\} \\
S^c := \{1,2,\cdots,n\} \setminus S.
\end{cases}
\end{equation}
Then, the equation (\ref{kernel}) is exactly the same equation
as we had in Theorem 4.20 in
\cite{cho:hsf03} with $(h_{i_s} -h_0)$ replaced by $c_{i_s}$.

Hence, by applying the same method, one can show that such
elements in the kernel of $\widetilde{\delta}_2$ lies in the image of
$\widetilde{\delta}_2$. This finishes the proof.
\end{proof}

Since for a fiber to have a non-trivial Floer cohomology is a very
special geometric property, it seems to deserves a name to them.

\begin{definition} We call {\it balanced} a Lagrangian fiber that
have a non-vanishing Floer cohomology.
\end{definition}

In the next section, we will provide a geometric description of
balanced torus fibers.

\section{Description of the balanced torus fibers}
\label{sec:fibers}

In this section, we now examine the equation (\ref{holboundary})
in terms of toric geometry. In particular, in the case of no line
bundle twisted, we provide a concrete toric description of the
conditions for a fiber to satisfy the equation.

For given $A \in int P$, we partition $G=G(\Sigma) = \{v_j \}_{1
\leq j \leq N}$ into the disjoint union
$$
G = \coprod_{\mu}G_{(A;\mu)}
$$
where $G_{(A;\mu)}$ is the set of $v_j \in G$ with the symplectic
area of the associated homotopy class $\beta_j = [D(v_j)] \in
\pi_2(X,L)$
$$
\omega_P(\beta_j) = \mu
$$
for each given positive number $\mu$. Obviously $G_{(A;\mu)} =
\emptyset$ except for a finite number of values of $\mu$'s
$$
0 < \mu_1 < \mu_2 < \cdots < \mu_{L_A}
$$
and $ 1 \leq L_A \leq N$. Then (\ref{holboundary}) becomes
\begin{equation}\label{holboundary2}
\sum_{v_j \in G_{(A;\mu_\ell)}} h^{v_j} v_j = \sum_{v_j \in
G_{(A;\mu_\ell)}} e^{i\langle \nu, v_j \rangle} v_j = 0
\end{equation}
for all $1 \leq \ell \leq L_A$.

\begin{prop} Assume $X_{\Sigma(P)}$ is Fano and
let $L = \mu^{-1}(A) \subset (X_{\Sigma(P)}, \omega_P)$ be a fiber
for $A \in \mbox{int }P$ and $\CL$ be a flat line bundle with the
holonomy vector $\nu = (\nu_1, \cdots, \nu_n)$ such that $A$ and
$\nu$ satisfy (\ref{holboundary2}). Then we have the isomorphism
$$
HF^{BM}((L,\CL);J_0) \cong H^*(L;\Lambda_{nov}^\CC).
$$
For all other cases, $HF^{BM}((L,\CL);J_0)$ is trivial.
\end{prop}

Now we specialize to the case without $\CL$, i.e., all
$h^{v_j}\equiv 1$. In the remaining section, we will provide a
more concrete description of the balanced fibers by analyzing
(\ref{holboundary2}) in terms of toric data.

Note that in this case (\ref{holboundary2}) just becomes
\begin{equation}\label{eq:sumvj}
\sum_{v_j \in G_{(A;\mu_\ell)}} v_j  = 0.
\end{equation}
We denote by
$$
\{1, \cdots, N\} = \coprod_{\ell=1}^{L_A} I_\ell
$$
the partition of $\{1, \cdots, N\}$ corresponding to the partition
of $G = \coprod_{\ell =1}^{L_A} G_{(A;\mu_\ell)}$. We also denote
\begin{equation}
e_{I_\ell}  =  \sum_{j \in I_\ell} e_j.
\end{equation}
By the exact sequence
$$
0 \to \KK \stackrel{i}\to \ZZ^N \stackrel{\pi}\to \ZZ^n \to 0
$$
(\ref{eq:sumvj}) implies that there exists $\delta_\ell \in \KK$
such that
$$
i(\delta_\ell) = e_{I_\ell} \in \ZZ^{I_\ell}
$$
for each $1 \leq \ell \leq L_A$, where $\ZZ^{I_\ell}$ is the
obvious product space. We denote by
$$
\Delta_\ell \subset (S^1)^{I_\ell}
$$
the obvious diagonal circle group generated by the vector
$e_{I_\ell}\in \ZZ^{I_\ell}$ and by $\Delta$ their products as a
subgroup of $(S^1)^N$. By construction, we have
$$
\Delta \subset K
$$
with $\dim \Delta = L_A \leq \dim K = N-n$.

We will now carry out the ``reduction by stages'' to describe our
toric manifolds $X_{\Sigma(P)}$ and the Lagrangian torus fiber $L$
in a two-step process. We denote $\delta = Lie(\Delta)$ and by
$$
\mu_\Delta: \CC^N \to \delta^*
$$
the moment map of the action of $\Delta$ on $\CC^N$. We denote
$$
j: \Delta \hookrightarrow K (\subset (S^1)^N) \quad\mbox{or }\, j:
\delta \hookrightarrow k (\subset \RR^N)
$$
the inclusion homomorphism and $T^\Delta:=(S^1)^N/\Delta$.

We note that $(S^1)^N$ acts on $\CC^N$ as the direct product of
the actions of $d_{\ell}$-dimensional torus $(S^1)^{I_\ell}$ on
$\CC^{I_\ell}$. By carrying out the {\it first} reduction by the
action of $\Delta$, we have obtained the reduced space
\begin{eqnarray*}
Y_\Delta & = &\mu_\Delta^{-1}(j^*(r))/\Delta \cong \PP^{(d_1-1)}
\times \cdots \times \PP^{(d_{L_A}-1)}\\
\omega_\Delta & = &\omega_1 \oplus \cdots \oplus \omega_{L_A}
\end{eqnarray*}
where we have
$$
j^*(r) = \sum_{\ell=1}^{L_A}(- \lambda_{I_\ell})e_{I_\ell}^* \in
\delta^*, \quad \lambda_{I_\ell} = \sum_{i\in I_\ell} (\lambda_i)
$$
with respect to the basis $\{e_{I_1}^*, \cdots, e_{I_\ell}^*\}$
dual to the basis $\{e_{I_1}, \cdots, e_{I_\ell}\}$ of $\delta$,
and $\omega_\ell$ is the Fubini-Study form on $\PP^{(d_\ell-1)}$
associated to the value $\lambda_{I_\ell}$ of the momentum
function $\mu_{\Delta_\ell}: \CC^{d_\ell} \to \delta_\ell^*\cong
\RR$, which becomes nothing but the standard momentum function of
the $S^1$ action on $\CC^{d_\ell}$ i.e.,
$$
z \in \CC^{d_\ell} \mapsto \frac{1}{2}|z|^2 \in \RR.
$$

Furthermore the residual torus $T^\Delta= (S^1)^N/\Delta$ is the
direct product
$$
T^\Delta = \prod_{\ell} T^{\Delta_\ell}
$$
with $T^{\Delta_\ell} := (S^1)^{I_\ell}/\Delta_\ell$, and
canonically acts on the reduced space $Y_\Delta$ as the direct
product action of the standard torus action of $T^{\Delta_\ell}
\cong (S^1)^{d_\ell}/\Delta_\ell$ on $\PP^{(d_\ell -1)}$. We
denote
$$
t^{\Delta_\ell} = Lie ((S^1)^{I_\ell}/\Delta_\ell)\cong
\RR^{(d_\ell-1)}.
$$
This action of the torus $T^{\Delta_\ell}$ on $\PP^{(d_\ell -1)}$
naturally extends to the action of the product
$$
U^\Delta : = \prod_{\ell}U(d_\ell)
$$
as the K\"ahler isometry with respect to the canonical complex and
symplectic structures induced from the ones on $Y_\Delta =
\mu_\Delta^{-1}(j^*(r))/\Delta$.

Now the quotient group $K/\Delta:= K_\Delta$ acts on $Y_\Delta$.
We denote its moment map by
$$
\mu_{K_\Delta}: Y_\Delta \to k_\Delta^*
$$
and the natural projection $K \to K_\Delta$ by $\pi_\Delta$. Then
we have the identity
$$
\pi_\Delta^*\circ \mu_{K_\Delta} = \mu_K.
$$
and the {\it second} reduction provides the description of
$(X_{\Sigma(P)}, \omega_P)$ as the reduced space
$$
X_{\Sigma(P)} \cong \mu_{K_\Delta}^{-1}(s)/K_\Delta
$$
where $s \in k_\Delta^*$ such that $\pi_\Delta^*(s) = r$.

In terms of this identification, the Lagrangian torus $L =
\mu^{-1}(A), \, A=(a_1, \cdots, a_n)$ can be written as
$$
\mu_{T^\Delta}^{-1}(A^1, \cdots, A^{L_A})/K_\Delta \cong
\Big(\mu_{T^{\Delta_1}}^{-1}(A^1) \times \cdots \times
\mu_{T^{\Delta_{L_A}}}^{-1}(A^{L_A})\Big)/K_\Delta
$$
where $A^\ell \in(t^{\Delta_\ell})^*$ and $(A^1, \cdots, A^{L_A})
\in \oplus_{\ell}(t^{\Delta_\ell})^*$ and
$$
\mu_{T^{\Delta_\ell}}: \PP^{(d_\ell -1)} \to (t^{\Delta_\ell})^*
\cong \RR^{(d_\ell -1)}
$$
is the standard moment map on $\PP^{(d_\ell-1)}$ of the action by
the torus $T^{\Delta_\ell}$. Here $(A^1, \cdots, A^{L_A}) =
\pi^*(A)$ where
$$
\pi: T^\Delta = (S^1)^N/\Delta \to \frac{(S^1)^N/\Delta}{K/\Delta}
\cong (S^1)^N/K =T^n.
$$
By the symmetry consideration, it follows that
$\mu_{T^{\Delta_\ell}}^{-1}(A^\ell)$ is the Clifford torus of
$\PP^{(d_\ell -1)}$.

We summarize the above discussion into the following theorem

\begin{theorem} Let $X_{\Sigma(P)} = \mu^{-1}(r)/K$ be a Fano
toric manifold with the canonical symplectic form $\omega_P$. Then
each balanced Lagrangian torus fiber in $X_{\Sigma(P)}$ has the
form
$$
L \cong (L_1 \times\cdots \times L_{L_A})/ K_\Delta \subset
\mu_{K_\Delta}^{-1}(s)/K_\Delta \cong X_{\Sigma(P)}
$$
where $K_\Delta = K/\Delta$ and $L_\ell$ is the Clifford torus of
$(\PP^{(d_\ell-1)}, \omega_\ell)$ with $\omega_\ell$ the
Fubini-Study form associated to the normalization
$$
\PP^{d_\ell-1} = \mu_{\Delta_\ell}^{-1}(-\lambda_{I_\ell})/S^1,
\quad \lambda_{I_\ell} = \sum_{i \in I_\ell} \lambda_i.
$$
\end{theorem}

\section{Hori-Vafa's B-Model Calculation}
\label{sec:hori}

In this section and the next, we will relate the equation
(\ref{holboundary}) with the critical point equation of the
superpotential of the Landau-Ginzburg mirror to the toric manifold
$(X_{\Sigma(P)},\omega_P)$, after substituting $T^{2\pi} =
e^{-1}$. We will closely follow the notations from [HV] with few
minor exceptions, and exclusively use convention that the letter
$i$ runs over $1, \cdots, N$, $a$ over $1, \cdots, k(=N-n)$ and
$\alpha$ over $1, \cdots, n$.

In this section, we first describe the prediction of Floer
cohomology by Hori via the mirror symmetry correspondence from
Hori and Vafa~\cite{hori:ms02} or Hori~\cite{hori:lmsd00}.

Suppose the $k$-dimensional torus $K =(S^1)^k$ acts on $\CC^N$ as
follows.
$$
(e^{i\theta_1},\cdots,e^{i\theta_k}) \cdot (z_1,\cdots,z_N) =
( \sum_{a=1}^k e^{iQ_{1a}\theta_a}z_1,\cdots, \sum_{a=1}^k
e^{iQ_{Na}\theta_a}z_N)
$$
The moment map of this action is given by
$$\mu : \CC^N \to (\RR^k)^* $$
$$ (z_1,\cdots,z_N) \mapsto \frac{1}{2}(\sum_i
Q_{i1}|z_i|^2,\cdots,\sum_i Q_{ik}|z_i|^2)
$$
(See section \ref{sec:forms} for more detailed discussion on
this). Now we consider the quotient $\mu^{-1}(r)/K$ as toric
manifolds where $r = (r_1, \cdots, r_k)$ lies in $(\RR^k)^*$.

With some physical arguments, Hori and Vafa [HV] introduce the
dual geometry by introducing periodic variables $Y_i,
i=1,\cdots,N$ with $Y_i \equiv Y_i + 2\pi i$ such that for
$a=1,\cdots,k$,
\begin{equation}\label{constraint}
\sum_{i=1}^{N} Q_{ia} Y_i = t_a
\end{equation}
where $t_a = r_a - i\theta_a$.

\begin{remark}
Here we consider the case where the B-field is zero
\end{remark}
The real part of $Y_i$ represents the position of the Lagrangian
torus fiber and imaginary part represents the holonomy of the line
bundle on this torus fiber. And one considers the superpotential
\begin{equation}\label{superpotential}
W:= \sum_{i=1}^N e^{-Y_i}.
\end{equation}
The critical points of the superpotential correspond to specific
fibers and holonomies whose Floer cohomology are non-vanishing.\\

For a given $Q$, we consider the equation
\begin{equation}\label{viQja}
\sum_{i=1}^N v_i Q_{ia} =0, \quad  v_i \in \ZZ.
\end{equation}
The space of solutions of (\ref{viQja}) form an integral lattice
of rank $n = N - k$ in $\RR^N$. We denote a $\ZZ$-basis of this
lattice by $\{v^\alpha\}_{1 \leq \alpha \leq n} \subset \RR^N$
with
$$
v^\alpha = (v^\alpha_1, \cdots, v^\alpha_N)
$$
each of them satisfying
\begin{equation}\label{qv}
\sum_{i=1}^N v_i^\alpha Q_{ia} =0.
\end{equation}
Therefore the general solutions for the constraint equation
$\sum_i Q_{ia} Y_i = t_a$ have the form
\begin{equation}\label{Y_i}
Y_i = \sum_{\alpha=1}^n v_i^\alpha \Theta_\alpha + y_i
\end{equation}
with $n=N-k$ periodic variables $\Theta_\alpha$ ($\mod 2\pi i$)
where $y = (y_1, \cdots, y_N)$ is a special solution of
$$\sum_{i=1}^N Q_{ia} Y_i = t_a.$$
(In [HV], the letters $t_i$'s are used for $y_i$'s which is
somewhat confusing with the other usage of $t_a$'s.)

Now the superpotential (\ref{superpotential}) of the mirror theory
can be expressed as
\begin{equation}\label{W}
W = \sum_{i=1}^N \exp (-y_i - \langle \Theta,v_i \rangle ),
\end{equation}
where
$$
v_i := (v_i^1, v_i^2, \cdots, v_i^n) \in \RR^n \cong N_\RR
$$
and $\langle \Theta,v_i \rangle $ is the short hand notation for
$\sum_{\alpha=1}^n v_i^\alpha \Theta_\alpha$. Note that the
condition $\frac{\partial W}{\partial \Theta_\alpha} =0$ is the
same as
\begin{equation}\label{criticaleq}
\sum_{i=1}^N e^{-Y_i}\cdot v_i^\alpha = \sum_{i=1}^N \exp (-y_i -
\langle \Theta,v_i \rangle ) \cdot v_i^\alpha =0
\end{equation}
for $\alpha = 1, \cdots, n$. One can already see the similarity
between equation (\ref{holboundary}) and the equation
(\ref{criticaleq}). In the next section, we show that two
equations indeed coincide, if we substitute
$$
T^{2\pi} = e^{-1}, \quad\mbox{and then}\quad  y_i = - \lambda_i.
$$

\section{Equivalence when $T^{2\pi}=e^{-1}$}
\label{sec:equiv}
In this section, we show that our calculation of
the (Bott-Morse) Floer cohomology indeed verifies the mirror
symmetry prediction made by Hori-Vafa's $B$-model calculation.
More precisely, the condition $(\ref{holboundary})$ to have
non-vanishing Floer cohomology with $T^{2\pi}=e^{-1}$ exactly
corresponds to the critical points of the superpotential, with a
canonical definition of the variables $Y_i$'s.

To see the correspondence for the compact toric manifold
$X_{\Sigma(P)}$, we define $Y_i$ as follows:

\begin{definition}
For $i=1,\cdots,N$, define $Y_i \in \RR \times i (\RR / 2\pi \ZZ)$
as
\begin{equation}
\begin{cases}
Re(Y_i) =  Area(\beta_i)/2\pi \\
Im(Y_i) = i\log( h^{v_i}) = -\langle \nu, v_i \rangle
\;\;\textrm{mod}\; 2\pi
\end{cases}
\end{equation}
where $v_i$'s are the generators of the one dimensional cones of
the fan $\Sigma$ associated the toric manifold $X_{\Sigma(P)}$
as in section \ref{sec:toric} and
$\beta_i \in \pi_2(X,L)$ is its associated homology class,
and $\nu$ is the holonomy vector of the
flat line bundle $\CL$ defined in (\ref{eq:nu}).
\end{definition}

Then it follows from Theorem \ref{area} that
\begin{equation}\label{eq:Yi}
Y_i =  (\langle A, v_i \rangle - \lambda_i) - i \langle \nu, v_i
\rangle = \langle A - i \nu, v_i \rangle - \lambda_i
\end{equation}
and hence
\begin{eqnarray}
e^{-Y_i} & = & e^{- (\langle A - i \nu, v_i \rangle - \lambda_i)}
\quad\mbox{or}
\label{eq:eYi1} \\
& = & h^{v_i} e^{-(\langle A, v_i \rangle - \lambda_i)}.
\label{eq:eYi2}
\end{eqnarray}

Consider the choice of $t_a$'s given by the real numbers
$$
t_a = - \sum_{i=1}^N Q_{ia}\lambda_i
$$
for $a =1,\cdots,N-n$. Then by the choice of $t_a$'s,
$$
y_i = - \lambda_i, \quad i = 1, \cdots, N
$$
is a special solution of (\ref{constraint}).

\begin{prop} For any vectors $A$ and $\nu$ above,
$Y_i$'s defined by (\ref{eq:Yi}) satisfy the constraint equation
$$
\sum_{i=1}^N Q_{ia} Y_i = t_a \;\;\mbox{for each}\;a = 1, \cdots,
N-n.
$$
\end{prop}
\begin{proof}
First note that we have the following equality from the exact
sequence (\ref{kexact}) or equation (\ref{qv})
$$
\sum_{i=1}^N Q_{ia} v_i^\alpha =0 \quad \mbox{ for all } \,\,
\alpha,
$$
and so we have
\begin{equation}
\sum_{i=1}^N Q_{ia} v_i =0.
\end{equation}
From this, we derive
\begin{eqnarray*}
\sum_{i=1}^N Q_{ia} Y_i
& =& \sum_{i=1}^N Q_{ia} (\langle A - i\nu, v_i \rangle -\lambda_i) \\
&=&  \langle A - i\nu, \sum_{i=1}^N Q_{ia}
v_i\rangle - \sum_{i=1}^N  Q_{ia}\lambda_i \\
&=& 0 + t_a = t_a
\end{eqnarray*}
which finishes the proof.
\end{proof}
Now identifying the variable $\Theta$
$$
\Theta = A - i\nu,
$$
$Y_i$'s defined in (\ref{eq:Yi}) coincide with (\ref{Y_i}).

Now, it remains to show that the condition (\ref{holboundary}) to
have non-vanishing Floer cohomology corresponds to the critical
points of the superpotential $W = \sum_{i=1}^N e^{-Y_i}$, if we
substitute $T^{2\pi}=e^{-1}$.
\begin{prop}
The $\Theta = A - i\nu$ is a critical point of the superpotential
$W$ if and only if $A$ and $\nu$ (or $h^{v_i}$'s) satisfy
(\ref{holboundary}), i.e., $\delta_2\langle pt \rangle = 0$.
\end{prop}
\begin{proof}
The condition (\ref{holboundary})
$$\sum_{i} h^{v_i} T^{2\pi(\langle A, v_i \rangle -\lambda_i)}\cdot v_i =0$$
becomes the following equation, after we substitute $T^{2\pi} =
e^{-1}$:
\begin{equation}\label{cri}
\sum_{i} h^{v_i} e^{-(\langle A, v_i \rangle -\lambda_i)}\cdot v_i
=0.
\end{equation}
Then from (\ref{eq:eYi2}) and from the choice $y_j = -\lambda_j$,
the above equation is same as
\begin{equation}
\sum_{i} e^{-Y_i}v_i^\alpha=0 \;\;\forall \; \alpha,
\end{equation}
which is precisely the condition for $Y_i$ to be the critical
points of the superpotential $W$ as in the equation
(\ref{criticaleq}).
\end{proof}

\section{Examples}

\subsection{The complex projective space $\PP^2$}
This example is taken from [Cho].  We consider $\PP^2$ associated
with the moment polytope $P$ defined by
\begin{equation}
\begin{cases}
\langle x,(1,0) \rangle \geq 0 \\
\langle x,(0,1) \rangle \geq 0 \\
\langle x,(-1,-1) \rangle \geq r
\end{cases}
\end{equation}
Let $(a_1,a_2) \in \mbox{int}(P)$. For the Lagrangian submanifold
$L:= \mu_T^{-1}(a_1,a_2)$, there exist three Maslov index 2 discs
(up to $Aut(D^2)$) with boundary in $L$.
 It is not hard to check that
its moment map $\mu_T$ trajectories are in fact straight lines. To
find the torus fiber whose Floer cohomology is non-vanishing, we
check the condition (\ref{holboundary}).
\begin{equation}
\delta_2(\langle pt \rangle)=h_1 T^{2\pi a_1} (1,0) + h_2 T^{2\pi
a_2} (0,1) + h_1^{-1} h_2^{-1} T^{2\pi (-a_1-a_2+r)} (-1,-1)=0
\end{equation}
Since $h_i \in U(1)$, we have
\begin{equation}
\begin{cases}
h_1 = h_1^{-1}h_2^{-2} \\
a_1 = -a_1-a_2+r
\end{cases}
\end{equation}
\begin{equation}
\begin{cases}
h_2 = h_1^{-1}h_2^{-2} \\
a_2 = -a_1-a_2+r
\end{cases}
\end{equation}
Hence, we have $$a_1 = a_2 = r/3,$$
$$h_1=h_2 \;\;\textrm{and} \;\; h_1^3=1$$
The Lagrangian fiber $\mu_T^{-1}(r/3,r/3)$ is called the
the {\em Clifford torus} and
the holonomies $(h_1,h_2)$ of the line bundle $\CL$ for the
non-vanishing Floer cohomologies on the Clifford torus are
$$ (1,1), (e^{2\pi i /3},e^{2\pi i /3}), (e^{4\pi i /3},e^{4\pi i
/3})$$
\subsection{Hirzebruch surfaces $F_c$}
This example illustrates well the difference between the actual
Floer cohomology with $\Lambda_{nov}$-coefficient and the Floer
cohomology with the parameter value $T^{2\pi} =e^{-1}$. The latter
was predicted by Hori-Vafa \cite{hori:ms02} by the B-model
calculation. However the latter version of the Floer cohomology is
{\it not} invariant under the Hamiltonian isotopy of the
Lagrangian torus fiber while the former version is so. The latter
Floer cohomology has Euler number of fibers whose Floer cohomology
is non-vanishing for toric Fano manifolds (\cite{hori:ms02}),
especially there exists four such fibers for Hirzebruch surfaces
$F_1$ and $F_2$. But we will show that the former version has {\it
no} fiber whose Floer cohomology is non-vanishing.

We start with the example $F_1$.
We consider $F_1$ associated with the moment polytope $P_1$
defined by
\begin{equation}
\begin{cases}
\langle x,(1,0) \rangle \geq -1 \\
\langle x,(0,1) \rangle \geq -1 \\
\langle x,(0,-1) \rangle \geq -1 \\
\langle x,(-1,1) \rangle \geq -1 \\
\end{cases}
\end{equation}
This polytope is {\em reflexive} \cite{B2}, thus
$F_1$ is Fano.
 Let
$(a_1,a_2) \in \mbox{int}(P_1)$. For the Lagrangian submanifold
$L:= \mu_T^{-1}(a_1,a_2)$, there exist four Maslov index 2 discs
(up to $Aut(D^2)$) with boundary in $L$. In the torus coordinates,
these are given as
$$(c_1,c_2z),(c_1z,c_2),(c_1,\frac{c_2}{z}),(\frac{c_1}{z},c_2z)$$
To find the torus fiber whose Floer cohomology is non-vanishing,
we check the condition (\ref{holboundary}).
\begin{eqnarray*}
0=\delta_2(\langle pt\rangle)=h_1 T^{2\pi (a_1+1)} (1,0) + h_2
T^{2\pi (a_2+1)}
(0,1) + \\
h_2^{-1} T^{2\pi (-a_2+1)} (0,-1) + h_1^{-1} h_2 T^{2\pi
(-a_1+a_2+1)} (-1,1)
\end{eqnarray*}
From the first coordinate, we obtain,
\begin{equation}
\begin{cases}
h_1 = h_1^{-1}h_2 \\
T^{2\pi(a_1+1)} = T^{2\pi(-a_1 +a_2 +1)}
\end{cases}
\end{equation}
Therefore we have $$h_2 = h_1^2, a_2 = 2a_1.$$ From the second
coordinate, we have
$$h_1^2 T^{2\pi(2a_1+1)} - h_1^{-2} T^{2\pi(-2a_1+1)} + h_1
T^{2\pi(a_1+1)}=0.$$ Or, equivalently
\begin{equation}\label{formaleq}
(h_1T^{2\pi a_1})^4 -1 + (h_1T^{2\pi a_1})^3 =0.
\end{equation}
Now, we substitute $T^{2\pi} = e^{-1}$. Then, the equation
(\ref{formaleq}) becomes
\begin{equation}
(h_1 e^{- a_1})^4 -1 + (h_1e^{-a_1})^3 =0.
\end{equation}
By setting $X:=h_1e^{- a_1} $, we have
\begin{equation}
X^4 + X^3 -1 =0.
\end{equation}
It is not hard to check the four solutions of this equation indeed
gives the location of the four fibers inside the polytope $P_1$,
whose Floer cohomology with the value $T^{2\pi} = e^{-1}$ is
non-vanishing. This agrees with the B-model calculation from
\cite{hori:ms02}.

For the Floer cohomology with $\Lambda_{nov}$-coefficient, note that
we regard $T$ as a formal parameter.
Hence to have a solution of the equation (\ref{formaleq}),
we should have
\begin{equation}
a_1 =0.
\end{equation}
In this cases, the equation becomes,
\begin{equation}\label{non-unitary}
h_1^4 + h_1^3 -1 =0.
\end{equation}
It is easy to check that this equation does not have a solution
for $h_1 \in U(1)$. Hence, there exists no torus fiber in
$F_1$(from the polytope $P_1$) whose Floer cohomology with
$\Lambda_{nov}$-coefficient is non-vanishing.

\begin{theorem}
Let $c$ be any positive integer and consider the Hirzebruch
surface $F_c$ with the canonical symplectic (K\"ahler)-form as a
toric manifold. Then all the non-singular torus fiber has trivial
Floer cohomology with $\Lambda_{nov}$-coefficients.
\end{theorem}
\begin{remark}
For $F_0 \cong \PP^1 \times \PP^1$ with both factors having the
same area, we find one fiber with four possible holonomies whose
Floer cohomology with $\Lambda_{nov}$-coefficient is
non-vanishing.
\end{remark}
\begin{proof}
We may consider a polytope $P_c$ (trapezoid with lengths $B,A,B+cA$) given as
follows:
\begin{equation}
\begin{cases}
\langle x,(1,0) \rangle \geq 0 \\
\langle x,(0,1) \rangle \geq 0 \\
\langle x,(0,-1) \rangle \geq -A \\
\langle x,(-1,-c) \rangle \geq -B-cA \\
\end{cases}
\end{equation}
Let $(a_1,a_2) \in \mbox{int}(P_c)$. Consider the Lagrangian
submanifold $L:= \mu_T^{-1}(a_1,a_2)$. We check the condition
(\ref{holboundary}) as before.
\begin{eqnarray*}
0=\delta_2(\langle pt\rangle)=h_1 T^{2\pi (a_1)} (1,0) + h_2
T^{2\pi (a_2)}
(0,1) + \\
h_2^{-1} T^{2\pi (-a_2+A)} (0,-1) + h_1^{-1} h_2^{-c} T^{2\pi
(-a_1-ca_2+B+cA)} (-1,-c)
\end{eqnarray*}
From the first coordinate, we have
\begin{equation}
\begin{cases}
h_1 = h_1^{-1}h_2^{-c} \\
a_1 = -a_1-ca_2+B+cA
\end{cases}
\end{equation}
And for $T$ as a formal parameter, from the second coordinate of
the above equation $\delta_2 \langle pt\rangle =0$, we have
\begin{equation}
a_2 = -a_2+A = -a_1-ca_2+B+cA.
\end{equation}
Combining these equations, we have
\begin{equation}
a_1=a_2=A/2,\;\; B = (\frac{2-c}{2})A.
\end{equation}
Hence, for $c \geq 2$, the length $B$ becomes non-positive which is not possible.
For the case $c=1$, we should have
$$h_1^4 +h_1^3 -1 =0,$$
which does not admit any solution in $U(1)$, in particular $h_1 =
1$ is not a solution. Hence, $\delta_2( \langle pt \rangle ) \neq
0$ and this finishes the proof.

\end{proof}

\begin{remark} The non-unitary solutions of (\ref{non-unitary})
for the $\Lambda_{nov}$-coefficients can be interpreted as the
solutions when a $B$-field is turned on: consider the complex
symplectic form
$$
\omega + \sqrt{-1} B=: \Omega
$$
and the one-parameter family $(X, z \Omega)$ with $z \in \CC^*$.
However the solutions for  $B\neq 0$ do not seem to carry any
natural symplectic geometrical meaning.
\end{remark}

\section{Obstruction classes, open Gromov-Witten invariants and
superpotentials} \label{sec:obstruct}

Fukaya-Oh-Ohta-Ono \cite{FOOO} have defined the obstruction cycles
of the filtered $A_\infty$-algebra associated to each Lagrangian
submanifold and developed a deformation theory thereof, which
tells whether one can kill the $\frak m_0$-term by a suitable
gauge equivalence. The $\frak m_0$ is defined by a collection of
currents induced by the (co)chains
$$
[\CM_1(\beta), ev_0]
$$
for all $\beta \in \pi_2(M,L)$. More precisely, we have
\begin{equation}\label{eq:m0}
\frak m_0(1) = \sum_{\beta \in \pi_2(M,L)}[\CM_1(\beta), ev_0]
\cdot T^{Area(\beta)} q^{\mu(\beta)/2} \in C^*(L) \otimes
\Lambda_{nov,0}.
\end{equation}
The sequence $\{o_k(L)\}_{1 \leq k <\infty}$ of the obstruction
classes introduced in [FOOO] is the iterative obstructions to
deforming the filtered $A_\infty$-structure so that
$$
\frak m_0 \equiv 0 \quad \mod T^{\lambda_{k+1}} \quad \mbox{as }
\, k \to \infty.
$$
Here we order those $\lambda$'s that appear as the area of
$\beta$, i.e., as $\omega(\beta)$.

Since the paper [FOOO] appeared, it became a folklore among some
mathematicians and physicists alike that under the mirror symmetry
correspondence FOOO's obstruction (co)chain in the $A$-model
should correspond to the superpotential in the $B$-model.

In fact, our computation confirms this test in the toric case. We
now explain this correspondence precisely. We first recall from
Theorem \ref{classify} and the orientability of the torus that
there is no holomorphic discs of Maslov index less than $2$.
According to [section 7, FOOO],  all obstruction classes
$o(\beta)$ are well-defined and the only non-trivial obstruction
classes (as currents) are the ones given by
\begin{equation}
o(\beta):=[\CM_1(\beta), ev_0] \quad \mbox{for $\beta$ with
$\mu(\beta) = 2$}
\end{equation}
for the torus fibers in this paper (see [section 7, FOOO] for more
explanation), which also coincides with $\frak m_0(1)$ in this
case. In view of Proposition \ref{disccord} and consideration of
the sign from [Cho], we also have
$$
o(\beta) = [L]\, (=1)
$$
the fundamental class of $L$ for any $\beta$ with $\mu(\beta)=2$.
Therefore we have obtained the formula for the obstruction class
of $L$
\begin{equation}\label{eq:o}
o(L) = \sum_{i=1}^N h^{v_j}T^{Area(\beta_j)}\cdot q
\end{equation}
from the definition of obstruction classes [Definition 4.6 \& 4.8,
FOOO]. However it follows, by the same substitution $T^{2\pi}$ by
$e^{-1}$ as before, that the right hand side of (\ref{eq:o})
precisely becomes
$$
\sum_{i=1}^N \exp(-y_i - \langle \Theta, v_i \rangle) = W(\Theta)
$$
if we ignore the harmless grading parameter $q$. Therefore we have
confirmed the exact correspondence
$$
o(L) \longleftrightarrow W
$$
for the case of Lagrangian torus fibers in toric manifolds.

In addition, comparing (\ref{holboundary}) with the derivative
$$
\frac{\partial W}{\partial \Theta}=\Big(\frac{\partial W}{\partial
\Theta_1}, \cdots,\frac{\partial W}{\partial \Theta_n}\Big)
$$
after substitution of $T^{2\pi} = e^{-1}$, we have also verified
the correspondence
$$
\delta_2\langle pt \rangle  \longleftrightarrow \frac{\partial
W}{\partial \Theta}.
$$
Recall from [addenda, O1] that the second named author observed
that the obstruction class is determined by the (genus zero) {\it
one-point} open Gromov-Witten invariants in the monotone case,
which obviously generalizes to the toric Fano case (because the
fiber does not have holomorphic discs of non-positive Maslov
indices). Similarly  in general, $\delta_2\langle pt \rangle $ is
determined by the {\it two-point} open Gromov-Witten invariants
from the definitions (\ref{del}) and (\ref{delall}) \cite{FOOO}.
Here we would like to emphasize that in this case of torus fibers
in the Fano toric manifolds, the (genus zero) open Gromov-Witten
invariants considered here are rigorously well-defined (with
respect to the canonical complex structure). This correspondence
between adding one marked point and taking the derivative is
consistent with the well-known principle in the calculus of
correlation functions in physics.

Combination of these facts suggests an intriguing relation between
the derivatives of $W$ of the superpotential (or of the
obstruction $o(L)$) and the ``open Gromov-Witten invariants'' of
$L$ in general. (Here we put the quotation mark because the open
Gromov-Witten invariants in general has not been rigorously
defined.) In fact, there has been conjectured by physicists that
the superpotential is related by the mirror symmetry
correspondence to the ``open Gromov-Witten potential'' (see
\cite{KKLM} for example), and our work provides a concrete
mathematical evidence via an $A$-model calculation. As far as we
understand, most calculations, if not all, in the physics
literature in this respect have been done in the $B$-model side.
We hope to further investigate this relation in the future.

\section{Discussion: non-Fano cases}
\label{sec:discussion}

We believe that our calculation of the Floer cohomology in this
paper remain to be true for the non-Fano toric manifolds. In this
section, we explain what remains to be proved for the non-Fano
cases.

The structure and regularity theorem of {\it smooth} holomorphic
discs still hold for the non-Fano case. However for {\it singular}
curves, distinction occurs in the trasversality problem because of
the presence of multiple covered spheres of {\it negative} Chern
numbers. Therefore it is essential to use the abstract
perturbation in the framework of Kuranishi structure \cite{FOn},
even if all disc components are already regular. With this
transversality problem taken care of, all the theorems, especially
those in sections \ref{sec:floer} and \ref{sec:computation} remain
to be true, possibly except the statement
$$
HF(L, \phi(L);J') \cong HF^{BM}(L,J_0).
$$
The proof of this isomorphism is expected to use a singular
degeneration argument as those used in \cite{FOh1} {\it in the
presence of non-trivial instantons} which may not be transversal.
Unlike the Fano case, there is no soft argument as those used in
[Oh4] to go from the limit configurations to the case of small
parameters because of non-transversality of sphere components. One
really has to construct a Kuranishi structure in the limit
configurations and to prove other non-trivial convergence
statements. This singular degeneration problem is currently being
studied by the second named author with K. Fukaya \cite{FOh2}.

In the end, we expect that the above isomorphism still holds but
details of the proof remain to be worked out. Because of this, we
restrict ourselves to the Fano case for sections \ref{sec:floer},
\ref{sec:computation} and the beginning of section
\ref{sec:index-two} in this paper.

\bibliographystyle{amsalpha}

\end{document}